\documentclass[11pt]{article}
\usepackage[a4paper,margin=28mm]{geometry}
\usepackage{amsmath,amssymb,amsthm,mathtools}
\usepackage{needspace}
\usepackage[hidelinks]{hyperref}
\numberwithin{equation}{section}
\newtheorem{theorem}{Theorem}[section]
\newtheorem{lemma}[theorem]{Lemma}
\newtheorem{proposition}[theorem]{Proposition}

\newcommand{\R}{\mathbb{R}}

\newcommand{\N}{\mathbb{N}}
\newcommand{\dd}{\,\mathrm{d}}
\newcommand{\norm}[1]{\lVert #1\rVert}
\newcommand{\abs}[1]{\lvert #1\rvert}

\newcommand{\dist}{\operatorname{dist}}

\hypersetup{pdftitle={The sharp Fourier extension estimate for the sphere},
pdfauthor={Pang-hung Chung}}
\title{The sharp Fourier extension estimate for the sphere}
\author{Pang-hung Chung}
\date{}
\begin{document}
\maketitle
\begin{abstract}
The full strong Fourier extension estimate for the unit sphere is established
in every ambient dimension at least three. The finite Knapp boundary is
obtained by combining a uniform diagonal estimate on rescaled spherical
graphs with a transverse product bound. A Whitney decomposition and a
dyadic dual pairing yield strong type at the exact scaling exponents.
The critical local estimate and the strict mixed range come from the
recent work on Euclidean local smoothing.
\end{abstract}

\section{Introduction}
\label{sec:intro}

Curvature enters the sphere restriction problem in two ways. For a
smooth density on the whole sphere, it gives decay of the Fourier
transform at large distances. For a density concentrated on a small
cap, it determines the shape of the region on which the oscillatory
integral remains coherent. The first effect controls the smallest
possible output exponent. The second relates the input exponent to the
output exponent.

For the unit sphere
\[
 S^{n-1}=\{\omega\in\R^n:|\omega|=1\}
\]
with Euclidean surface measure \(\sigma\), define
\[
 E_n f(X)=\int_{S^{n-1}}e^{2\pi i X\cdot\omega} f(\omega)\dd\sigma(\omega).
\]
The sphere restriction conjecture asks for the exact exponents for which
$E_n$ is bounded from $L^p(S^{n-1},\sigma)$ to $L^q(\R^n)$.
Since $\sigma$ is finite, $E_n f$ is continuous for every
$f\in L^p(\sigma)$ with $p\geq1$. We write $p'$ for the H\"older conjugate
of $p$, including $1'=\infty$ and $\infty'=1$.
For a measure space \((\Omega,\nu)\), a measurable \(f\), and
\(1\le r\le\infty\), put
\[
 \norm f_{L^r(\Omega,\nu)}
 =\begin{cases}
   \bigl(\int_\Omega |f|^r\dd\nu\bigr)^{1/r},&1\le r<\infty,\\
   \operatorname*{ess\,sup}_{\Omega,\nu}|f|,&r=\infty.
  \end{cases}
\]
The space \(L^r(\Omega,\nu)\) consists of finite-norm functions modulo
\(\nu\)-almost-everywhere equality.

The decay and cap conditions together characterize the full strong-type
range for the sphere.

\begin{theorem}\label{thm:main}
Let \(n\ge3\) and \(1\le p,q\le\infty\). There exists
\(C_{n,p,q}<\infty\) such that
\[
 \norm{E_n f}_{L^q(\R^n)}
 \le C_{n,p,q}\norm{f}_{L^p(S^{n-1},\sigma)}
 \qquad \bigl(f\in L^p(S^{n-1},\sigma)\bigr)
\]
if and only if
\begin{equation}\label{eq:region}
 q>\frac{2n}{n-1}
 \qquad\text{and}\qquad
 q\geq\frac{n+1}{n-1}\,p'.
\end{equation}
\end{theorem}

To see the distinction between the two conditions in
\eqref{eq:region}, first take \(f=1\). On a sequence of radial intervals of fixed length,
the extension has modulus comparable to
\(\abs X^{-(n-1)/2}\). At \(q=2n/(n-1)\), integrating its \(q\)-th
power over these intervals gives a divergent harmonic sum.
Now let \(f\) be the indicator of a spherical cap of radius \(h\).
The cap has measure comparable to \(h^{n-1}\), and its extension has
modulus comparable to \(h^{n-1}\) on a region with \(n-1\) tangential
dimensions of size \(h^{-1}\) and one normal dimension of size
\(h^{-2}\). Thus the cap test gives
\begin{equation}\label{eq:cap-scaling-intro}
 \frac{\norm{E_n f}_{L^q(\R^n)}}{\norm{f}_{L^p(S^{n-1},\sigma)}}
 \gtrsim h^{(n-1)(1-1/p)-(n+1)/q}.
\end{equation}
The exponent in \eqref{eq:cap-scaling-intro} must be nonnegative as
\(h\) tends to zero. At equality the cap test gives a bound of order
one from below, leaving open the question of sufficiency.

Multilinear restriction estimates provide another way to exploit
curvature and transversality. Bennett, Carbery and Tao established the
multilinear restriction and Kakeya estimates in 2006 \cite{BCT}.
Bourgain and Guth used these estimates in 2011 to obtain bounds for
linear restriction and oscillatory integral operators \cite{BG}.
Polynomial partitioning subsequently gave a restriction estimate for
positively curved surfaces in three dimensions \cite{guth2016}, followed
by higher-dimensional estimates based on broad norms \cite{guth2018}.
These developments explain the role of separated cap interactions in
the restriction problem. In the boundary argument below, the required
uniformity is obtained directly for rescaled spherical phases.

The Tomas--Stein theorem and interpolation with \(L^1\to L^\infty\)
give the part of this line with \(p\leq2\) \cite{stein,tomas}.
For positively curved surfaces in \(\R^3\), Kim obtained the portion
\(q>13/4\) \cite[Theorem~1.2]{kim}. The Euclidean local smoothing
results of OpenAI give a critical local bound with arbitrarily small
power loss and global strong type when both inequalities in
\eqref{eq:region} are strict
\cite[Lemma~9.2 and Corollary~9.3]{openai-local-smoothing}.
A limit of the strict-range estimates would require uniform control of
their constants as the exponents approach the line. Moreover, at the
line the scaling power in \eqref{eq:cap-scaling-intro} vanishes, so
single-cap scaling provides no geometric decay for the sum over cap
sizes. The following boundary estimate supplies the remaining
sufficiency in Theorem~\ref{thm:main}.

\begin{theorem}[Finite Knapp boundary]\label{thm:boundary}
Let $n\geq3$ and $1\leq p,q\leq\infty$. If
\begin{equation}\label{eq:boundary}
 \frac{2n}{n-1}<q=\frac{n+1}{n-1}\,p'<\infty,
\end{equation}
then \(\exists\,C_{n,p,q}<\infty\) with
\[
 \norm{E_n f}_{L^q(\R^n)}\leq C_{n,p,q}\norm{f}_{L^p(S^{n-1},\sigma)}
 \qquad (f\in L^p(S^{n-1},\sigma)).
\]
\end{theorem}

To prove the boundary estimate, we separate interactions between caps
according to their relative scale, following the bilinear argument of
Tao, Vargas and Vega \cite{TVV}. Tao
obtained sharp estimates for separated paraboloid patches
\cite{tao-paraboloid}. Buschenhenke describes the passage to the
Knapp line by interpolation between diagonal and bilinear estimates
\cite[Section~1.1]{buschenhenke}. In our setting the estimates used
after rescaling must hold uniformly for spherical graph patches
whose sizes and centres both vary.

Write \(d=n-1\) and \(s_c=2(d+1)/d\). At the local exponent \(s_c\),
the physical Jacobian and the input scaling combine with the
normalisation of the integral to give the cap power
\[
 -d+\frac{d+2}{s_c}+\frac{d}{s_c}=0.
\]
The local estimate still has a loss because the physical observation
ball grows under rescaling. For caps small relative to the observation
scale, the phase is close to a quadratic on that ball. The comparison
with the paraboloid in \cite[Section~12]{openai-local-smoothing}
illustrates this mechanism. Section~\ref{sec:analytic-inputs} compares
the spherical phases directly at auxiliary scales and then removes the
local loss with constants uniform in the scale and centre.

Uniform diagonal estimates control a single rescaled patch. To sum
interactions between patches, we use separation to gain room in the
input exponent. The \(L^2\) product estimate in
\cite[Theorem~2.3]{TVV} illustrates the required gain.
Section~\ref{sec:transverse} derives a uniform spherical estimate by
reducing the collision geometry to a slab bound. Interpolation with
the diagonal estimate carries the gain to the product exponent needed
for \eqref{eq:boundary}.

In the Whitney sum the relevant comparison is between the volume
\(v_j\) of a cell at scale \(j\) and the total mass \(M=\abs A\) of
the set supporting the restricted input. For \(v_j\ge M\), the
total mass bound gives a positive power of \(M/v_j\). For \(v_j<M\),
the bound by the volume of each individual cell gives a positive
power of \(v_j/M\). The two geometric sums give the exact restricted
estimate, as quantified in \eqref{kn:two-scale}. Finally, two nearby
restricted estimates on the Knapp line give decay on opposite sides
of a dyadic pairing, yielding the strong estimate at the desired
exponents.

\subsection*{Conventions and organization}

For \(k\ge1\), the Euclidean ball in \(\R^k\) and the orthogonal group are
\[
 B(a,r)=\{x\in\R^k:|x-a|<r\},\qquad
 O(k)=\{Q\in\R^{k\times k}:Q^{\mathsf T}Q=I_k\}.
\]
For functions on \(\R^k\), write
\[
 \nabla=(\partial_{x_1},\ldots,\partial_{x_k}),\qquad
 \Delta=\sum_{j=1}^k\partial_{x_j}^2.
\]
The notation \(D^j\phi\) denotes the \(j\)-th Fr\'echet derivative.
Write \(\N=\{0,1,2,\ldots\}\) and
\[
 \mathbf1_A(x)=
 \begin{cases}1,&x\in A,\\0,&x\notin A.\end{cases}
\]
Within a proof, \(\norm f_r\) abbreviates the \(L^r\) norm on the
domain and measure specified for \(f\) in that proof or its statement.
For subsets of one Euclidean space,
\[
 \dist(A,B)=\inf\{|a-b|:a\in A,\ b\in B\},\qquad \inf\varnothing=\infty.
\]
Functions under integrals are measurable on their indicated domains.
In graph coordinates Lebesgue measure is used, and \(|A|\) denotes
the measure of a measurable set \(A\). The notation $A\lesssim B$ means $A\leq C B$,
with the dependence of $C$ specified in the surrounding statement.
For nonnegative quantities,
\[
 a\asymp_{\mathcal P}b
 \quad\text{means}\quad
 a\lesssim_{\mathcal P}b\ \text{and}\ b\lesssim_{\mathcal P}a.
\]
The notation \(R=O_{\mathcal P}(b)\) means
\(\norm R\le C_{\mathcal P}b\) in the relevant Euclidean or tensor norm.
Constants denoted by $C$ or $C$ with subscripts may change between
estimates and depend only on the indicated parameters.
The graph phase is $x\cdot v+t\phi(v)$, and the spherical extension
uses the kernel \(e^{2\pi iX\cdot\omega}\) in the definition of \(E_n\).
We establish necessity and the
uniform graph estimates first. The transverse estimate and Whitney
sum yield the strong estimate on a graph patch. A finite cap covering
then gives Theorem~\ref{thm:boundary}, from which we obtain
Theorem~\ref{thm:main}.

\section{Necessary conditions}
\label{sec:necessity}

Fix \(n\ge3\), put \(d=n-1\), and write \(E=E_n\) in this section. The radial condition is strict because the oscillatory leading term of
\(E1\) is nonzero on intervals of fixed length at arbitrarily large radii.
The other condition is determined by coherent propagation from a shrinking
cap. We record both arguments with the normalization of Theorem~\ref{thm:main}.
For the constant density, the radial Helmholtz equation identifies
the leading oscillation.

\begin{lemma}[Nonzero radial oscillation]
\label{nec:constant-asymptotic}
Let \(d=n-1\). There exist \(A,B\in\R\) such that
\[
 A^2+B^2>0,\qquad
 E1\left(\frac{r}{2\pi}e\right)
 =r^{-d/2}\bigl(A\cos r+B\sin r+O_n(r^{-1})\bigr)
\]
for all \(e\in S^{n-1}\) and \(r\ge1\).
\end{lemma}

\begin{proof}
Fix \(e\in S^{n-1}\). Put
\[
  U(x)=\int_{S^d}e^{ix\cdot\omega}\,\dd\sigma(\omega),
  \qquad F(r)=U(re).
\]
Rotation and reflection invariance give
\[
 U(Qx)=U(x)\quad(Q\in O(n)),\qquad \overline{U(x)}=U(-x)=U(x).
\]
Differentiation under the integral yields
\[
 \Delta U(x)
 =-\int_{S^d}|\omega|^2e^{ix\cdot\omega}\,\dd\sigma(\omega)
 =-U(x).
\]
The radial Laplacian therefore gives
\begin{equation}\label{nec:radial-ode}
  F''(r)+\frac{d}{r}F'(r)+F(r)=0,\qquad r>0 .
\end{equation}
Set
\[
 \beta=\frac d2,\qquad u(r)=r^\beta F(r),\qquad
 V(r)=\frac{\beta(\beta-1)}{r^2}.
\]
Equation~\eqref{nec:radial-ode} becomes
\begin{equation}\label{nec:conjugated-ode}
 u''+u=V(r)u.
\end{equation}
The energy \(H=u^2+(u')^2\) satisfies
\begin{equation}\label{nec:energy-differential}
  H'=2Vuu',\qquad -|V|H\le H'\le |V|H .
\end{equation}
Also \(H(1)>0\), since
\(F(1)=\int_{S^d}\cos(e\cdot\omega)\,\dd\sigma(\omega)>0\).
Integrating the inequalities in \eqref{nec:energy-differential} yields
\[
  H(1)\exp\left(-\int_1^r |V(s)|\,\dd s\right)
  \le H(r)
  \le H(1)\exp\left(\int_1^r |V(s)|\,\dd s\right).
\]
Since \(\int_1^\infty |V(s)|\,\dd s<\infty\),
\begin{equation}\label{nec:energy-bounds}
 0<\inf_{r\ge1}H(r)\le \sup_{r\ge1}H(r)<\infty.
\end{equation}
Set
\[
  a=u\cos r-u'\sin r,
  \qquad b=u\sin r+u'\cos r .
\]
Then \(a^2+b^2=H\), and \eqref{nec:conjugated-ode} gives
\begin{equation}\label{nec:coefficient-derivatives}
  a'=-Vu\sin r,\qquad b'=Vu\cos r .
\end{equation}
Equations~\eqref{nec:coefficient-derivatives} and~\eqref{nec:energy-bounds} imply \(a',b'\in L^1([1,\infty))\), with limits
\(A=\lim_{r\to\infty}a(r)\) and \(B=\lim_{r\to\infty}b(r)\). The lower bound in \eqref{nec:energy-bounds} gives
\[
 A^2+B^2=\lim_{r\to\infty}H(r)\ge\inf_{r\ge1}H(r)>0.
\]
The derivative identities \eqref{nec:coefficient-derivatives} give
\begin{equation}\label{nec:coefficient-tail}
 \begin{aligned}
 |a(r)-A|+|b(r)-B|
 &\le\int_r^\infty\bigl(|a'(s)|+|b'(s)|\bigr)\,\dd s\\
 &\lesssim_n\int_r^\infty s^{-2}\,\dd s=r^{-1}.
 \end{aligned}
\end{equation}
Finally,
\[
 E1\left(\frac r{2\pi}e\right)
 =F(r)=r^{-d/2}\bigl(a(r)\cos r+b(r)\sin r\bigr),
\]
and \eqref{nec:coefficient-tail} gives the asserted expansion.
\end{proof}

The nonzero amplitude gives the strict radial condition. A shrinking
cap supplies the second necessary condition.

\begin{proposition}
\label{nec:necessary}
Let \(1\le p,q\le\infty\). If
\(E:L^p(S^{n-1})\to L^q(\R^n)\) is bounded, then
\[
 \begin{cases}
 q>2n/(n-1),\\
 (d+2)/q\le d(1-1/p).
 \end{cases}
\]
For every \(1\le p\le\infty\),
\[
 \norm{Ef}_{L^\infty(\R^n)}
 \le \sigma(S^{n-1})^{1-1/p}\norm f_{L^p(S^{n-1})}.
\]
\end{proposition}

\begin{proof}
\emph{The radial condition.}
The constant function belongs to every \(L^p(S^d)\).
Let \(L(r)=A\cos r+B\sin r\) be the leading sinusoid of
Lemma~\ref{nec:constant-asymptotic}. Choose \(c>0\) and an interval
\((s,s+w)\) of length \(0<w<2\pi\) on which \(|L(r)|\ge 2c\).
For all sufficiently large \(k\) and
\[
  I_k=(s+2\pi k,s+w+2\pi k),
\]
the asymptotic gives
\(|E1(re/(2\pi))|\ge c r^{-d/2}\) on \(I_k\).
For finite \(q\), polar integration and disjointness give
\begin{equation}\label{nec:radial-integrability}
 \norm{E1}_{L^q(\R^n)}^q
 \gtrsim_{n,q}\sum_{k\gg1}\int_{I_k}r^{d-dq/2}\,\dd r.
\end{equation}
For \(q\le2n/d\), the exponent satisfies \(d-dq/2\ge-1\), and hence
\[
 \sum_{k\gg1}\int_{I_k}r^{d-dq/2}\,\dd r
 \gtrsim_{n,q}\sum_{k\gg1}k^{d-dq/2}
 \ge\sum_{k\gg1}k^{-1}=\infty.
\]
The divergence in \eqref{nec:radial-integrability} forces \(q>2n/d\). This inequality also holds for
\(q=\infty\).

\smallskip
\emph{The Knapp condition.}
For \(0<\delta\le 1/2\) and \(v\in B^d(0,1)\), write
\[
  h_\delta(v)=(1-\delta^2|v|^2)^{1/2},
  \qquad \omega_\delta(v)=(\delta v,h_\delta(v)).
\]
The area element is \(\delta^d h_\delta(v)^{-1}\,\dd v\). Define
\[
 f_\delta(\omega_\delta(v))=h_\delta(v),\qquad
 f_\delta=0\quad\text{off }\omega_\delta(B^d(0,1)).
\]
For finite \(p\),
\[
  \norm{f_\delta}_{L^p(S^d)}^p
  =\delta^d\int_{B^d(0,1)} h_\delta(v)^{p-1}\,\dd v
  \le |B^d(0,1)|\,\delta^d .
\]
For \(p=\infty\), \(\norm{f_\delta}_{L^\infty}\le1\). Thus
\begin{equation}\label{nec:cap-input}
 \norm{f_\delta}_{L^p(S^d)}\lesssim_{n,p}\delta^{d/p}
 \qquad(1\le p\le\infty).
\end{equation}
Use the physical coordinates
\[
  T_\delta(x,t)=\left(\frac{x}{2\pi\delta},
                      -\frac{t}{\pi\delta^2}\right),
  \qquad
  \phi_\delta(v)=\frac{2(1-h_\delta(v))}{\delta^2}
                =\frac{2|v|^2}{1+h_\delta(v)} .
\]
Cancellation of the area density by \(f_\delta\) gives
\begin{equation}\label{nec:cap-phase}
  Ef_\delta(T_\delta(x,t))=
  \delta^d e^{-2it/\delta^2}
  \int_{B^d(0,1)} e^{i(x\cdot v+t\phi_\delta(v))}\,\dd v .
\end{equation}
On \(U=B^d(0,1/8)\times(-1/8,1/8)\),
\[
 |\phi_\delta(v)|\le2,\qquad
 |x\cdot v+t\phi_\delta(v)|\le\tfrac18+\tfrac14=\tfrac38 .
\]
Writing \(B=B^d(0,1)\), the inequality \(|e^{is}-1|\le|s|\) gives
\[
 \left|\int_B e^{i(x\cdot v+t\phi_\delta(v))}\,\dd v-|B|\right|
 \le\tfrac38|B|.
\]
Hence
\begin{equation}\label{nec:cap-coherence}
 \left|\int_B e^{i(x\cdot v+t\phi_\delta(v))}\,\dd v\right|
 \ge\tfrac58|B|\ge\tfrac12|B|.
\end{equation}
The set \(T_\delta U\) has volume
\begin{equation}\label{nec:cap-volume}
  |T_\delta U|=\frac{2}{(2\pi)^{d+1}}\delta^{-(d+2)}|U|.
\end{equation}
Equations~\eqref{nec:cap-phase}, \eqref{nec:cap-coherence}
and~\eqref{nec:cap-volume} give the lower bound below.
Boundedness and \eqref{nec:cap-input} give the upper bound:
\begin{equation}
\label{nec:knapp}
  \delta^{d-(d+2)/q}
  \lesssim_{n,p,q} \norm{Ef_\delta}_{L^q(\R^n)}
  \lesssim_{n,p,q} \delta^{d/p} .
\end{equation}
For \(q=\infty\), the lower bound uses the essential supremum on
\(T_\delta U\), with \(1/\infty=0\). Equation~\eqref{nec:knapp} gives
\[
 \delta^{d(1-1/p)-(d+2)/q}\lesssim_{n,p,q}1
 \quad(0<\delta\le\tfrac12),
\]
so
\[
 d(1-1/p)-(d+2)/q\ge0.
\]
In particular \(p=1\) forces \(q=\infty\).

\smallskip
\emph{The infinite-output estimate.}
The modulus of the kernel is one, so H\"older's inequality on the
finite-measure sphere yields
\begin{equation}
\label{nec:infty}
  \norm{Ef}_{L^\infty(\R^n)}
  \le \norm{f}_{L^1(S^d)}
  \le \sigma(S^d)^{1-1/p}\norm{f}_{L^p(S^d)} .
\end{equation}
\end{proof}

The infinite-output case is thus settled. For finite output exponents,
we begin with the critical local estimate and its global consequence
in the strict range.

\section{Local restriction and uniform estimates for spherical graphs}
\label{sec:analytic-inputs}

Throughout this section, let \(d=n-1\), and put
\[
 b_c=\frac{2(d+1)}{d+2},\qquad s_c=b_c'=\frac{2(d+1)}d.
\]
We use the critical local restriction estimate from
\cite[Lemma~9.2]{openai-local-smoothing}. In the sphere normalization of
this paper, its dual form is
\begin{equation}\label{eq:critical-local-sphere}
 \norm{E_n g}_{L^{s_c}(B(0,R))}
 \le C_{n,\epsilon}R^\epsilon\norm{g}_{L^{s_c}(S^{n-1})},
 \qquad R\ge1,\quad \epsilon>0.
\end{equation}

\Needspace{7\baselineskip}
Epsilon removal and interpolation give the following global estimate
when both exponent inequalities are strict.

\begin{theorem}[Strict mixed range, \cite{openai-local-smoothing}]
\label{input:strict}
Let
\[
 1<p\le\infty,\qquad q<\infty,\qquad
 q>s_c,\qquad q>(s_c-1)p'.
\]
There exists \(C_{n,p,q}<\infty\) such that
\[
 \norm{E_n g}_{L^q(\R^n)}
 \le C_{n,p,q}\norm{g}_{L^p(S^{n-1})}
 \qquad(g\in L^p(S^{n-1})).
\]
\end{theorem}

\begin{proof}
This is the sphere case of \cite[Corollary~9.3]{openai-local-smoothing}.
Epsilon removal gives
\begin{equation}\label{eq:strict-diagonal}
 \norm{E_n g}_{L^s(\R^n)}
 \lesssim_{n,s}\norm g_{L^s(S^{n-1})}
 \qquad(s_c<s<\infty).
\end{equation}
Choose \(s\) and define \(r\) by
\[
 s_c<s<\min\{q,1+q/p'\},\qquad
 r=\begin{cases}
 (s-1)p',&p<s,\\
 s,&p\ge s.
 \end{cases}
\]
Interpolating \eqref{eq:strict-diagonal} with \eqref{nec:infty}
when \(p<s\), and using finite surface measure when \(p\ge s\), gives
\begin{equation}\label{eq:strict-mixed-intermediate}
 \norm{E_n g}_{L^r(\R^n)}
 \lesssim_{n,p,s}\norm g_{L^p(S^{n-1})},
 \qquad r<q.
\end{equation}
With \(\theta=r/q\), equations \eqref{eq:strict-mixed-intermediate}
and \eqref{nec:infty} give
\[
 \norm{E_n g}_{L^q}
 \le\norm{E_n g}_{L^r}^{\theta}
       \norm{E_n g}_{L^\infty}^{1-\theta}
 \lesssim_{n,p,q}\norm g_{L^p(S^{n-1})}.
\]
\end{proof}

The boundary argument also needs a diagonal estimate that is uniform
under changes in the size and centre of a graph patch. Let \(B=B^d(0,1)\).
For a real phase \(\phi\) on \(B\), write
\[
 \mathcal E_\phi g(x,t)=\int_B e^{i(x\cdot u+t\phi(u))}g(u)\dd u.
\]
For \(0<\delta,r\le1/2\), \(\abs{a}\le1/2\), define
\begin{equation}\label{eq:spherical-graph-phase}
 \Phi_\delta(v)=\frac{2}{\delta^2}
       \bigl(1-\sqrt{1-\delta^2\abs{v}^2}\bigr).
\end{equation}
The phase rebased at \(a\) and scale \(r\) is
\begin{equation}\label{eq:rebased-graph-phase}
 \Psi_{\delta,a,r}(u)
 =\frac{\Phi_\delta(a+ru)-\Phi_\delta(a)
         -r\nabla\Phi_\delta(a)\cdot u}{r^2}.
\end{equation}

The phases \(\Psi_{\delta,a,r}\) describe the patches produced by the
rescalings below. The next bound is uniform over this family.

\begin{theorem}[Uniform diagonal estimate]\label{input:diagonal}
For \(s_c<s<\infty\), there exists \(C_{n,s}<\infty\) such that
\[
 \norm{\mathcal E_{\Psi_{\delta,a,r}}g}_{L^s(\R^{d+1})}
 \le C_{n,s}\norm{g}_{L^s(B)}
 \qquad(g\in L^s(B)),
\]
uniformly under
\[
 0<\delta,r\le\tfrac12,\qquad |a|\le\tfrac12.
\]
\end{theorem}

The proof of Corollary~12.1 in \cite{openai-local-smoothing} compares
a spherical cap of radius chosen in terms of \(R\) with a fixed
paraboloid on \(B(0,R)\). Here the Whitney rescalings require control
for every \((\delta,a,r)\), with a constant independent of the triple.
Indeed, writing \(h=\delta r\), a direct change of variables in the
global spherical diagonal bound at \(s>s_c\) incurs the factor
\(h^{-d+(2d+2)/s}\). Its exponent is negative.

We work first at \(s_c\), where this exponent is zero. The change
of variables in the local estimate leaves a loss of order
\(R^\eta h^{-2\eta}\). Fix \(M>1+(d+1)/s_c\). At scales
\(h\ge R^{-M}\), choosing \(\eta\) small absorbs this loss in
\(R^\epsilon\). For smaller scales the rebased phase differs from its
quadratic part by \(O_n(\delta^2 r)\). The local error from this phase
difference is bounded by \(C_n R^{1+(d+1)/s_c}h\norm{g}_{s_c}\).
The quadratic part reduces to the paraboloid by a uniformly bounded
linear change of variables, followed by a fixed dilation of the
parameter ball. The local estimate (12.4) of
\cite{openai-local-smoothing} therefore also gives the bound at small
scales. The proof below makes the comparison separately in \(\delta\) and
\(r\), using auxiliary spherical phases at both steps. We then apply
the sparse decomposition method of \cite{tao-eps} directly to the
graph measure with density \(1_B\), retaining the same uniformity.

\subsection{Uniform local estimates}

At the critical exponent \(s_c\), the cap power cancels. The first
step is to control the remaining radius loss uniformly.

\begin{lemma}\label{lem:uniform-local-graphs}
Let
\[
 R\ge1,\qquad\epsilon>0,\qquad
 \phi\in\{\Phi_\delta,\Psi_{\delta,a,r}\},\qquad
 0<\delta,r\le\tfrac12,\quad |a|\le\tfrac12.
\]
Then
\[
 \norm{\mathcal E_\phi g}_{L^{s_c}(B^{d+1}(0,R))}
 \le C_{n,\epsilon}R^\epsilon\norm{g}_{L^{s_c}(B)}
 \qquad(g\in L^{s_c}(B)).
\]
\end{lemma}

\begin{proof}
The cap parameterization is
\[
 P_\delta(v)=(\delta v,\sqrt{1-\delta^2\abs v^2}),
\]
with surface density \(\delta^d(1-\delta^2\abs v^2)^{-1/2}\).
On the cap set
\[
 h(P_\delta(v))=g(v)\sqrt{1-\delta^2\abs v^2},
\]
and extend \(h\) by zero. At
\(X=(x/(2\pi\delta),-t/(\pi\delta^2))\), definition
\eqref{eq:spherical-graph-phase} gives
\begin{equation}\label{eq:cap-extension-identity}
 E_n h(X)=\delta^d e^{-2it/\delta^2}\mathcal E_{\Phi_\delta}g(x,t).
\end{equation}
For \(T_\delta(x,t)=X\), the change of variables satisfies
\begin{align*}
 |\det T_\delta|&\asymp_n\delta^{-(d+2)},&
 T_\delta(B(0,R))&\subset B(0,C_nR/\delta^2),\\
 \norm h_{s_c}&\lesssim_n\delta^{d/s_c}\norm g_{s_c}.
\end{align*}
The critical exponents satisfy
\begin{equation}\label{eq:critical-scaling-cancellation}
 -d+\frac{d+2}{s_c}+\frac d{s_c}=0.
\end{equation}
Combining \eqref{eq:cap-extension-identity} with
\eqref{eq:critical-local-sphere} gives
\[
 \norm{\mathcal E_{\Phi_\delta}g}_{L^{s_c}(B(0,R))}
 \lesssim_{n,\eta}
 \delta^{-d+(d+2)/s_c}(R/\delta^2)^\eta
 \delta^{d/s_c}\norm g_{s_c},
\]
By \eqref{eq:critical-scaling-cancellation}, this becomes
\begin{equation}\label{eq:cap-local-with-loss}
 \norm{\mathcal E_{\Phi_\delta}g}_{L^{s_c}(B(0,R))}
 \lesssim_{n,\eta} R^\eta\delta^{-2\eta}\norm{g}_{s_c}
 \qquad (\eta>0).
\end{equation}

For phases \(\phi,\widetilde\phi\) on \(B\), set
\(\kappa=1+(d+1)/s_c\). The pointwise inequality
\(\abs{e^{iv}-e^{iw}}\le\abs{v-w}\) and H\"older's inequality give
\begin{equation}\label{eq:local-phase-comparison}
 \norm{(\mathcal E_\phi-\mathcal E_{\widetilde\phi})g}_{L^{s_c}(B(0,R))}
 \lesssim_n R^\kappa\norm{\phi-\widetilde\phi}_{L^\infty(B)}
                 \norm{g}_{s_c}.
\end{equation}
Uniformly on \(B\), definition \eqref{eq:spherical-graph-phase} gives
\begin{equation}\label{eq:quadratic-cap-approximation}
 \norm{\Phi_\delta-|\,\cdot\,|^2}_{L^\infty(B)}
 \lesssim_n\delta^2.
\end{equation}
If \(\delta^2R^\kappa\ge1\), then
\[
 R^\eta\delta^{-2\eta}\le R^{(1+\kappa)\eta}.
\]
If \(\delta^2R^\kappa<1\), set
\(\delta_*=\tfrac12R^{-\kappa/2}\). Equation~\eqref{eq:local-phase-comparison}
and \eqref{eq:quadratic-cap-approximation} give
\begin{equation}\label{eq:auxiliary-cap-comparison}
 R^\kappa\norm{\Phi_\delta-\Phi_{\delta_*}}_\infty
 \lesssim_n R^\kappa(\delta^2+\delta_*^2)\lesssim_n1.
\end{equation}
For \(\delta^2R^\kappa\ge1\), use \eqref{eq:cap-local-with-loss}
at \(\delta\). For \(\delta^2R^\kappa<1\), use
\eqref{eq:cap-local-with-loss} at \(\delta_*\), followed by
\eqref{eq:local-phase-comparison} and \eqref{eq:auxiliary-cap-comparison}.
Both parameter ranges give
\begin{equation}\label{eq:uniform-cap-local}
 \norm{\mathcal E_{\Phi_\delta}g}_{L^{s_c}(B(0,R))}
 \lesssim_{n,\eta}\bigl(1+R^{(1+\kappa)\eta}\bigr)\norm g_{s_c}.
\end{equation}
Taking \(\eta=\epsilon/(1+\kappa)\) proves the first assertion.

For the rebased phase, define \(g_{a,r}\) by
\[
 g_{a,r}(a+ru)=g(u)\quad(u\in B),\qquad
 g_{a,r}=0\quad\text{on }B\setminus B(a,r),
\]
and set \(\ell=\nabla\Phi_\delta(a)\).
Definition \eqref{eq:rebased-graph-phase} gives
\begin{equation}\label{eq:rebased-extension-identity}
\begin{aligned}
 &\mathcal E_{\Phi_\delta}g_{a,r}
       (r^{-1}x-r^{-2}t\ell,r^{-2}t)\\
 &\qquad
 =r^d e^{i(r^{-1}x\cdot a+r^{-2}t(\Phi_\delta(a)-\ell\cdot a))}
       \mathcal E_{\Psi_{\delta,a,r}}g(x,t).
\end{aligned}
\end{equation}
The physical map \(T_{a,r}(x,t)=(r^{-1}x-r^{-2}t\ell,r^{-2}t)\)
satisfies
\begin{align*}
 |\det T_{a,r}|&=r^{-(d+2)},&
 T_{a,r}(B(0,R))&\subset B(0,C_nR/r^2),\\
 \norm{g_{a,r}}_{s_c}&=r^{d/s_c}\norm g_{s_c}.
\end{align*}
Apply \eqref{eq:uniform-cap-local}, with its loss exponent chosen as
\(\eta/(1+\kappa)\), on \(B(0,C_nR/r^2)\).
Equations \eqref{eq:rebased-extension-identity} and
\eqref{eq:critical-scaling-cancellation} then give
\begin{equation}\label{eq:rebased-local-with-loss}
 \norm{\mathcal E_{\Psi_{\delta,a,r}}g}_{L^{s_c}(B(0,R))}
 \lesssim_{n,\eta}R^\eta r^{-2\eta}\norm g_{s_c}.
\end{equation}
For \(H=\sqrt{1-\delta^2\abs v^2}\), direct differentiation gives
\[
 D^2\Phi_\delta(v)=2H^{-1}I+2\delta^2 H^{-3}vv^{\mathsf T},
 \qquad \norm{D^3\Phi_\delta(v)}\lesssim_n \delta^2
 \quad(\abs v\le1).
\]
Taylor's theorem and \eqref{eq:rebased-graph-phase} give
\begin{equation}\label{eq:rebased-quadratic-remainder}
 \left\|\Psi_{\delta,a,r}
 -\tfrac12 u^{\mathsf T}D^2\Phi_\delta(a)u\right\|_{L^\infty(B)}
 \lesssim_n\delta^2r.
\end{equation}
If \(rR^\kappa\ge1\), then
\[
 R^\eta r^{-2\eta}\le R^{(1+2\kappa)\eta}.
\]
For \(rR^\kappa<1\), set \(r_*=\tfrac12R^{-\kappa}\).
Subtracting the two quadratic approximations in
\eqref{eq:rebased-quadratic-remainder} gives
\begin{equation}\label{eq:auxiliary-rebased-comparison}
 R^\kappa\norm{\Psi_{\delta,a,r}-\Psi_{\delta,a,r_*}}_\infty
 \lesssim_n R^\kappa\delta^2(r+r_*)\lesssim_n1.
\end{equation}
For \(rR^\kappa\ge1\), apply \eqref{eq:rebased-local-with-loss}
at \(r\). For \(rR^\kappa<1\), combine
\eqref{eq:local-phase-comparison},
\eqref{eq:auxiliary-rebased-comparison} and
\eqref{eq:rebased-local-with-loss} at \(r_*\). This gives
\[
 \norm{\mathcal E_{\Psi_{\delta,a,r}}g}_{L^{s_c}(B(0,R))}
 \lesssim_{n,\eta}\bigl(1+R^{(1+2\kappa)\eta}\bigr)\norm g_{s_c}.
\]
Taking \(\eta=\epsilon/(1+2\kappa)\) completes the proof.
\end{proof}

\subsection{Restriction on sparse families}

To remove the radius loss, we pass to the dual restriction estimate
on well-separated physical balls.

For the remainder of the proof of Theorem~\ref{input:diagonal},
write \(m=d+1\), \(b=b_c\), \(\psi=\Psi_{\delta,a,r}\), and
\(\Gamma=\{(u,\psi(u)):u\in B\}\). Let \(\mu\) be the pushforward
of \(1_B\dd u\) to \(\Gamma\). The Fourier transform in this
subsection is
\[
 \widehat f(\xi)=\int_{\R^m}e^{-iz\cdot\xi}f(z)\dd z,\qquad
 \widehat\mu(z)=\int_\Gamma e^{-iz\cdot\omega}\dd\mu(\omega).
\]
All the measure estimates below are uniform in the graph parameters.

We have \(D^2\psi(u)=D^2\Phi_\delta(a+ru)\ge2I\).
The phase and its first two derivatives are uniformly bounded, so
the graph measures have common bounded support and bounded mass.
The graph measure satisfies
\begin{equation}\label{eq:graph-mass}
 \mu(B(\xi,h))\lesssim_n h^d.
\end{equation}
Its Fourier transform satisfies
\begin{equation}\label{eq:graph-mass-decay}
 \abs{\widehat\mu(z)}\lesssim_n (1+\abs z)^{-1/2}.
\end{equation}
Equation~\eqref{eq:graph-mass} follows by projection. For decay at \(z=(x,t)\),
slice \(B\) into intervals parallel to a coordinate axis.
If \(\abs x\le C_n\abs t\), the phase has second derivative of
magnitude at least \(2\abs t\) on each interval. The second-derivative
van der Corput estimate gives \(O_n(\abs t^{-1/2})\). If \(\abs x>C_n\abs t\), choose
\(j\) with \(\abs{x_j}\ge\abs x/\sqrt d\).
For sufficiently large \(C_n\), the first derivative in that coordinate
has magnitude at least \(\abs{x_j}/2\) and is monotone.
The first-derivative estimate gives \(O_n(\abs x^{-1})\) after
integrating over the transverse slice variables. Bounded \(z\) is
covered by the total mass bound.

Write \(T_h=\bigcup_{\omega\in\Gamma}B(\omega,h)\) for \(0<h\le1\).
A graph point within distance \(h\) of \(\xi\in T_h\) has a graph
neighborhood of radius \(h\) carrying measure \(\gtrsim_n h^d\).
Indeed, its projection contains the intersection of \(B\) with a
parameter ball of radius \(c_n h\) about the projected point.
It follows that \(\mu(B(\xi,2h))\gtrsim_n h^d\). Lemma~\ref{lem:uniform-local-graphs},
duality, modulation, and Fubini give, for \(f\) supported in \(B(0,R)\),
\begin{align}
 \int_{T_h}\abs{\widehat f(\xi)}^b\dd\xi
 &\lesssim_n h^{-d}
    \int_{\abs v<2h}\int_{\Gamma}
          \abs{\widehat f(\omega+v)}^b\dd\mu(\omega)\dd v
 \notag\\
 &\lesssim_{n,\eta} h R^{b\eta}\norm f_{L^b(\R^m)}^b
 \qquad (R\ge1,\ \eta>0).
 \label{eq:thin-neighborhood}
\end{align}

The thin-neighbourhood bound \eqref{eq:thin-neighborhood}
controls each ball separately. The
Fourier decay in \eqref{eq:graph-mass-decay} controls interactions
between balls and yields the following estimate.

\begin{lemma}\label{lem:graph-sparse}
Let \(N\in\N\) and \(z_1,\ldots,z_N\in\R^m\). Assume
\[
 D=16m,\qquad R\ge1,\qquad
 |z_i-z_j|\ge(RN)^D\quad(1\le i\ne j\le N).
\]
For every \(\eta>0\) and \(f\in L^b(\R^m)\) with
\[
 f=0\quad\text{a.e.\ on }\R^m\setminus\bigcup_{j=1}^N B(z_j,R),
\]
one has
\[
 \norm{\widehat f}_{L^b(\mu)}
 \lesssim_{n,\eta}R^\eta\norm f_{L^b(\R^m)}.
\]
\end{lemma}

\begin{proof}
For \(N=0\), the function \(f\) is zero almost everywhere.
Assume \(N\ge1\). Choose a Schwartz function \(w\) whose Fourier transform is
nonnegative, even, smooth and supported in \(B(0,1/8)\), and
such that \(w\) is bounded away from zero on \(B(0,1)\).
Set
\[
 f_j=f1_{B(z_j,R)},\qquad
 g_j(y)=\frac{f_j(z_j+y)}{w(y/R)}\quad(|y|<R),
\]
and extend \(g_j\) by zero. Disjointness and the lower bound for \(w\) give
\begin{equation}\label{eq:sparse-data-norm}
 \sum_j\norm{g_j}_b^b\lesssim\sum_j\norm{f_j}_b^b=\norm f_b^b.
\end{equation}
The Fourier product identity represents \(\widehat f|_\Gamma\),
up to a fixed normalization factor, as the image of the functions
\(\widehat g_j|_{T_{1/R}}\) under the operator with kernels
\[
 K_j(\omega,\xi)=e^{-iz_j\cdot\omega}R^m
                     \widehat w\bigl(R(\omega-\xi)\bigr).
\]
Let \(T\) be this kernel operator, with counting measure in \(j\),
Lebesgue measure in \(\xi\in T_{1/R}\), and \(\mu\) in \(\omega\).
Equation~\eqref{eq:graph-mass} gives
\begin{equation}\label{eq:sparse-kernel-mass}
 \sup_{j,\xi}\int_\Gamma |K_j(\omega,\xi)|\,\dd\mu(\omega)
 \lesssim_n R,\qquad
 \sup_{j,\omega}\int_{T_{1/R}}|K_j(\omega,\xi)|\,\dd\xi
 \lesssim_n1.
\end{equation}
Hence the diagonal Gram blocks have absolute row and column
integrals \(O_n(R)\).

For \(i\ne j\), define the off-diagonal Gram kernel by
\[
 G_{ij}(\xi,\zeta)=\int e^{i(z_i-z_j)\cdot\omega}
      R^{2m}\overline{\widehat w(R(\omega-\xi))}
                 \widehat w(R(\omega-\zeta))\dd\mu(\omega).
\]
The product of the two bumps is the transform of a convolution of
two modulated or reflected copies of \(w(\cdot/R)\).
Each has weighted \(L^1\) norm \(O(R^{m+1/2})\) with weight
\((1+\abs y)^{1/2}\), and that weight is submultiplicative.
Fourier inversion and \eqref{eq:graph-mass-decay} give
\begin{equation}\label{eq:sparse-gram-off-diagonal}
 |G_{ij}(\xi,\zeta)|
 \le C_nR^{2m+1}(1+\abs{z_i-z_j})^{-1/2}
 \lesssim_n N^{-1}.
\end{equation}
The set \(T_{1/R}\) has bounded volume.
Equations \eqref{eq:sparse-kernel-mass} and
\eqref{eq:sparse-gram-off-diagonal}, followed by the Schur test, give
\begin{equation}\label{eq:sparse-operator-endpoints}
 \begin{aligned}
 \norm T_{\ell^2_jL^2(T_{1/R})\to L^2(\mu)}&\lesssim_nR^{1/2},\\
 \norm T_{\ell^1_jL^1(T_{1/R})\to L^1(\mu)}&\lesssim_nR.
 \end{aligned}
\end{equation}
Interpolating \eqref{eq:sparse-operator-endpoints} gives
\begin{equation}\label{eq:sparse-operator-b}
 \norm T_{\ell^b_jL^b(T_{1/R})\to L^b(\mu)}
 \lesssim_n R^{1/b}.
\end{equation}
Combining \eqref{eq:sparse-operator-b},
\eqref{eq:thin-neighborhood} at \(h=1/R\), and
\eqref{eq:sparse-data-norm} gives
\begin{align*}
 \norm{\widehat f}_{L^b(\mu)}
 &\lesssim_nR^{1/b}
   \left(\sum_j\int_{T_{1/R}}|\widehat g_j(\xi)|^b\,\dd\xi\right)^{1/b}\\
 &\lesssim_{n,\eta}R^{1/b}R^{-1/b+\eta}
             \left(\sum_j\norm{g_j}_b^b\right)^{1/b}
 \lesssim_{n,\eta}R^\eta\norm f_b.
\end{align*}
This proves the claim.
\end{proof}

\subsection{Removal of the local loss}

Grouping unit balls into sparse families allows
Lemma~\ref{lem:graph-sparse} to yield an arbitrarily small power of
the number of balls.

\begin{proof}[Proof of Theorem~\ref{input:diagonal}]
First, if \(f\) is supported in a union of \(M\ge1\) unit balls, then
for every \(\epsilon>0\)
\begin{equation}\label{eq:unit-ball-loss}
 \norm{\widehat f}_{L^b(\mu)}
 \lesssim_{n,\epsilon} M^\epsilon\norm f_{L^b(\R^m)}.
\end{equation}
Let \(E\) be the set of centers. Choose a positive integer \(L\) and set
\[
 K=\lceil M^{1/L}\rceil,\qquad
 R_0=4(M+1),\qquad R_{j+1}=R_j^{3D+1}.
\]
Every center \(z\) belongs to at least one layer \(0\le j<L\)
characterized by
\begin{equation}\label{eq:sparse-cover-layers}
 \#(E\cap B(z,2R_j))\ge K^j,\qquad
 \#(E\cap B(z,2R_{j+1}))\le K^{j+1}.
\end{equation}
This follows by taking the first index \(j+1\) at which the second
inequality holds, using \(K^L\ge M\). Assign each center to one
such layer. In layer \(j\) choose a maximal \(4R_j\)-separated
subfamily. Its \(5R_j\)-balls cover the corresponding unit balls.
The selected centers within distance \(R_{j+1}\) of a fixed
selected center have disjoint \(2R_j\)-balls in the fixed center's
\(2R_{j+1}\)-ball. The two cardinality bounds in \eqref{eq:sparse-cover-layers}
bound their number by \(K\). Coloring the graph with these adjacencies takes at most
\(K+1\) colors. Each color gives a sparse family of \(5R_j\)-balls
since \(R_{j+1}\ge(5R_jM)^D\). A measurable partition subordinate to
the resulting cover and Lemma~\ref{lem:graph-sparse} bound the
restriction norm by
\begin{equation}\label{eq:sparse-cover-cost}
 \norm{\widehat f}_{L^b(\mu)}
 \le C_{n,L,\eta} M^{1/L+\eta(3D+1)^L}\norm f_b.
\end{equation}
Choose
\[
 L>\frac2\epsilon,\qquad
 0<\eta<\frac{\epsilon}{2(3D+1)^L}.
\]
Then \(1/L+\eta(3D+1)^L<\epsilon\), so
\eqref{eq:sparse-cover-cost} proves \eqref{eq:unit-ball-loss}.

If the \(z_j\) are unit separated and \(\abs{a_j}\le1\), apply
\eqref{eq:unit-ball-loss} to
\(\sum_j a_j1_{B(z_j,\rho)}\).
Choose a common small \(\rho\) so that the transform of
\(1_{B(0,\rho)}\) is bounded away from zero on the common measure
support and the translated balls are disjoint. Then
\begin{equation}\label{eq:bounded-phase-sums}
 \left\|\sum_{j=1}^M a_j e^{-iz_j\cdot\xi}\right\|_{L^b(\mu)}
 \lesssim_{n,\epsilon} M^{1/b+\epsilon}.
\end{equation}
For \(1<q<b\), normalize \(\sum_j|a_j|^q=1\) and set
\[
 J_k=\{j:2^{-k-1}<|a_j|\le2^{-k}\}\quad(k\ge0).
\]
Then
\begin{equation}\label{eq:dyadic-coefficient-mass}
 \#J_k\,2^{-q(k+1)}
 \le\sum_{j\in J_k}|a_j|^q\le1,\qquad
 \#J_k\le2^{q(k+1)}.
\end{equation}
Apply \eqref{eq:bounded-phase-sums} to \(2^ka_j\) on \(J_k\), and
use \eqref{eq:dyadic-coefficient-mass}. This gives
\begin{align*}
 \left\|\sum_j a_je^{-iz_j\cdot\xi}\right\|_{L^b(\mu)}
 &\le\sum_{k\ge0}2^{-k}
       \left\|\sum_{j\in J_k}2^ka_je^{-iz_j\cdot\xi}\right\|_{L^b(\mu)}\\
 &\lesssim_{n,\epsilon}
 \sum_{k\ge0}2^{-k(1-q(1/b+\epsilon))}<\infty
\end{align*}
when \(q(1/b+\epsilon)<1\). Finite measure and homogeneity give
\begin{equation}\label{eq:graph-phase-sums}
 \left\|\sum_j a_j e^{-iz_j\cdot\xi}\right\|_{L^q(\mu)}
 \lesssim_{n,q}\left(\sum_j\abs{a_j}^q\right)^{1/q}.
\end{equation}
For \(f\in C_c^\infty(\R^m)\), tiling gives
\begin{equation}\label{eq:lattice-fourier-tiling}
 \widehat f(\xi)
 =\int_{[0,1)^m}\sum_{k\in\mathbb Z^m}
       f(u+k)e^{-i(u+k)\cdot\xi}\,\dd u.
\end{equation}
The points \(u+k\) are unit separated. Apply Jensen's inequality to
\eqref{eq:lattice-fourier-tiling}, and then use
\eqref{eq:graph-phase-sums}. This gives
\[
 \norm{\widehat f}_{L^q(\mu)}^q
 \lesssim_{n,q}\int_{[0,1)^m}\sum_{k\in\mathbb Z^m}\abs{f(u+k)}^q\dd u
 = \norm f_{L^q(\R^m)}^q .
\]
Density extends this inequality to \(L^q(\R^m)\). With \(q=s'\),
\[
 s>s_c
 \quad\Longleftrightarrow\quad 1<s'<b_c,
\]
and duality gives
\[
 \norm{\mathcal E_\psi g}_{L^s(\R^m)}
 \lesssim_{n,s}\norm g_{L^s(B)}.
\]
This is the claimed diagonal estimate.
\end{proof}

The diagonal estimate controls each factor of an extension product.
We next use separation of the input supports to improve the mass
exponent for the product.

\section{A transverse gain for spherical graphs}
\label{sec:transverse}

Fix \(d\geq2\) and write \(D=d/(d+1)\).
We use \(\Phi_\delta\), \(\Psi_{\delta,a,r}\), and
\(\mathcal E_\phi\) from Section~\ref{sec:analytic-inputs}, with
\(0<\delta,r\leq1/2\) and \(|a|\leq1/2\).
We use
\[
 x_K(y)=1-\frac{d+2}{2d}y
\]
for the reciprocal input exponent on the Knapp line when the
product output exponent is \(1/y\). For their normalised elliptic
class, Tao, Vargas and Vega \cite[Theorem~2.3]{TVV} prove the
separated product estimate in \(L^2\) at reciprocal input
\(x_{\mathrm{TVV}}=3/4-1/(2(d+1))\). In these coordinates,
\[
 x_{\mathrm{TVV}}-x_K(1/2)=\frac{1}{2d(d+1)}>0 .
\]
Our spherical proof uses the same collision representation
\cite[Section~2.8]{TVV}. An exact stereographic identity expresses
the collision condition as a Euclidean slab condition and gives an
estimate uniform for \(0<\delta,r\leq1/2\), \(|a|\leq1/2\).
The interpolation
below requires only a positive gap beyond \(x_K(1/2)\) at \(L^2\).

\subsection{An estimate for slabs}

For \(h>0\) and a measurable \(V\subset\R^d\) of finite measure, define
\[
 S_h V(\omega,s)
   =h^{-1}\bigl|\{v\in V: |s-v\cdot\omega|\leq h\}\bigr|,
 \qquad \omega\in S^{d-1},\quad s\in\R.
\]
We use surface measure in \(\omega\) and Lebesgue measure in \(s\).

We need an \(L^m\) bound for these averages that is uniform in the
slab width.

\begin{lemma}\label{trans:slab}
Let
\[
 2\leq m<d+1,\qquad L>0,\qquad M_0<\infty.
\]
For every measurable \(V\subset\R^d\) with \(|V|\leq M_0\),
\[
 \norm{S_h V}_{L^m(S^{d-1}\times[-L,L])}
 \lesssim_{d,m,L,M_0}|V|^{1-1/m}
 \qquad(h>0).
\]
\end{lemma}

\begin{proof}
Set \(F=\mathbf1_V\), \(M=|V|\), and use the Fourier transform
\(\widehat F(\xi)=\int e^{-2\pi i v\cdot\xi}F(v)\dd v\).
The Gaussian kernel satisfies
\[
 G_h(s)=h^{-1}e^{-\pi(s/h)^2},\qquad
 \widehat G_h(\rho)=e^{-\pi h^2\rho^2},\qquad
 G_h(s)\geq e^{-\pi}h^{-1}\mathbf1_{\{|s|\leq h\}}.
\]
Hence
\begin{align}
 S_hV(\omega,s)
 &\leq e^\pi\int_{\R^d}F(v)G_h(s-v\cdot\omega)\dd v\notag\\
 &=e^\pi\int_{\R}e^{-\pi h^2\rho^2}
           \widehat F(\rho\omega)e^{2\pi i s\rho}\dd\rho.
 \label{trans:slab-fourier}
\end{align}
The integral on the right is real and nonnegative. The Fourier
identity follows by Fubini since \(F\in L^1\) and the Gaussian
multiplier is integrable.

Split the frequency integral in \eqref{trans:slab-fourier} into \(|\rho|\leq1\) and
\(2^j<|\rho|\leq2^{j+1}\), \(j\geq0\). Its low-frequency part
is bounded in absolute value by \(2M\). Define
\begin{align*}
 T_j(\omega,s)
 &=\int_{2^j<|\rho|\leq2^{j+1}}
       e^{-\pi h^2\rho^2}\widehat F(\rho\omega)e^{2\pi i s\rho}\dd\rho,\\
 \mathcal A_j&=\{\xi:2^j<|\xi|\leq2^{j+1}\}.
\end{align*}
Then
\begin{equation}\label{trans:slab-infty}
 \norm{T_j}_\infty\lesssim 2^jM.
\end{equation}
Plancherel in \(s\) and polar coordinates in \(\xi\) give
\begin{equation}\label{trans:slab-two}
 \norm{T_j}_{L^2(S^{d-1}\times\R)}^2
 \leq 2\cdot2^{-(d-1)j}
       \int_{\mathcal A_j}|\widehat F(\xi)|^2\dd\xi
 \leq2\cdot2^{-(d-1)j}M.
\end{equation}
The factor two in \eqref{trans:slab-two} accounts for both signs of
\(\rho\). Interpolating \eqref{trans:slab-infty} and
\eqref{trans:slab-two} gives
\begin{align*}
 \norm{T_j}_{L^m(S^{d-1}\times\R)}
 &\leq \norm{T_j}_\infty^{1-2/m}
        \norm{T_j}_{L^2(S^{d-1}\times\R)}^{2/m}\\
 &\lesssim (2^jM)^{1-2/m}
       \bigl(2^{-(d-1)j/2}M^{1/2}\bigr)^{2/m}\\
 &=2^{j(1-(d+1)/m)}M^{1-1/m}.
\end{align*}
The low-frequency term has norm \(O_{d,L,m}(M)\). Since
\[
 M\leq M_0^{1/m}M^{1-1/m},\qquad
 \sum_{j\geq0}2^{j(1-(d+1)/m)}
   =\frac{1}{1-2^{1-(d+1)/m}}<\infty,
\]
the finite frequency sums satisfy the claimed bound uniformly.
The Gaussian multiplier gives pointwise convergence, and Fatou's
lemma passes the bound to the limit.
\end{proof}

\subsection{The four-point identity}

We now identify the slabs that arise from the four-point phase defect.

Write \(\phi=\Psi_{\delta,a,r}\), \(k=\delta r\), and
\[
 b(v)=\sqrt{1-\delta^2|a+rv|^2},\qquad
 P(v)=(\delta(a+rv),b(v)).
\]
For \(|v|\leq1\) the point \(P(v)\) is on \(S^d\) and
\(3/4\leq b(v)\leq1\). The horizontal projection and the
identity
\[
 b(v)-b(w)
 =\frac{\delta^2\bigl(|a+rw|^2-|a+rv|^2\bigr)}{b(v)+b(w)}
\]
give the chord bounds
\begin{equation}\label{trans:chord-bounds}
 k|v-w|\leq |P(v)-P(w)|\leq2k|v-w|
 \qquad (|v|,|w|\leq1).
\end{equation}

Fix \(|u|\leq1\), put \(p=P(u)\), and define
\[
 p^\perp=\{Y\in\R^{d+1}:Y\cdot p=0\}.
\]
The stereographic coordinates in \(p^\perp\) are
\begin{equation}\label{trans:stereo-chart}
 \Xi_u(v)
   =k\left(\frac{P(v)-p}{|P(v)-p|^2}+\frac p2\right),
   \qquad v\ne u.
\end{equation}
Tangency follows from
\(p\cdot(P(v)-p)=-|P(v)-p|^2/2\). Euclidean inversion applied to
\eqref{trans:stereo-chart} gives
\begin{equation}\label{trans:stereo-distance}
 |\Xi_u(v)-\Xi_u(w)|
 =\frac{k|P(v)-P(w)|}{|P(v)-p|\,|P(w)-p|}.
\end{equation}
The same definition gives
\begin{equation}\label{trans:stereo-radius}
 |\Xi_u(v)|^2=\frac{k^2}{|P(v)-p|^2}-\frac{k^2}{4}.
\end{equation}
For \(\sigma>0\), put
\(V_{u,\sigma}=\{v:|v|\leq1,\ |v-u|\geq\sigma\}\).
Equations~\eqref{trans:chord-bounds} and \eqref{trans:stereo-distance}
give
\[
 \frac1{16}|v-w|
 \leq|\Xi_u(v)-\Xi_u(w)|
 \leq\frac2{\sigma^2}|v-w|
 \qquad(v,w\in V_{u,\sigma}).
\]
Equations~\eqref{trans:chord-bounds} and \eqref{trans:stereo-radius}
also give
\[
 \Xi_u(V_{u,\sigma})\subset\{Y:1/8\leq|Y|\leq1/\sigma\}.
\]
After an orthogonal identification \(p^\perp\simeq\R^d\), the
bi-Lipschitz bounds and the Hausdorff-measure characterization of
Lebesgue measure imply that there is \(C_{d,\sigma}\geq1\) such that
\begin{equation}\label{trans:stereo-measure}
 C_{d,\sigma}^{-1}|A|\leq|\Xi_u(A)|\leq C_{d,\sigma}|A|
 \qquad(A\subset V_{u,\sigma}\ \text{measurable}).
\end{equation}
The image \(\Xi_u(V_{u,\sigma})\) is compact. Borel sets map to
relative Borel sets, and the bi-Lipschitz bounds preserve null sets.
Thus \eqref{trans:stereo-measure} also applies to Lebesgue measurable
sets and gives the corresponding integral bounds for nonnegative
functions.

For \(z=u+v-w\) with \(|u|,|v|,|w|,|z|\leq1\), define
\begin{equation}\label{trans:phase-defect}
 \Delta(u,v,w)=\phi(u)+\phi(v)-\phi(w)-\phi(u+v-w).
\end{equation}
Put \(H=b(u)+b(v)-b(w)+b(z)\).
The constants and linear terms in \(\phi\) cancel, so
\[
 \Delta H
 =-\frac{2}{k^2}\bigl(|P(u)+P(v)-P(w)|^2-1\bigr).
\]
Here the horizontal coordinate of \(P(u)+P(v)-P(w)\) equals
that of \(P(z)\), and \(H\) is between \(1\) and \(3\).
If \(v,w\ne u\), let \(Y=\Xi_u(v)\), \(W=\Xi_u(w)\),
\(V_1=P(v)-p\), and \(W_1=P(w)-p\). Expanding the square and
the stereographic inner product gives
\begin{align*}
 |P(u)+P(v)-P(w)|^2-1
   &=2\bigl(|W_1|^2-V_1\cdot W_1\bigr),\\
 Y\cdot(W-Y)
   &=k^2\frac{V_1\cdot W_1-|W_1|^2}{|V_1|^2|W_1|^2}.
\end{align*}
Thus
\begin{equation}\label{trans:factor}
 \Delta H=4\left(\frac{|V_1|}{k}\right)^2
            \left(\frac{|W_1|}{k}\right)^2Y\cdot(W-Y).
\end{equation}
For \(|v-u|,|w-u|\geq\sigma\) and \(|\Delta|\leq\varepsilon\),
equations~\eqref{trans:chord-bounds} and \eqref{trans:factor} give,
with \(\omega=Y/|Y|\) and \(s=|Y|\),
\begin{equation}\label{trans:collision-slab}
 |s-W\cdot\omega|\leq6\varepsilon/\sigma^4.
\end{equation}

The slab estimate therefore bounds the measure of configurations
with small phase defect. We also record a Gaussian version for use
in the product estimate. Integrands containing \(\Delta\) below are
assigned value zero when \(|u+v-w|>1\), in accordance with their
displayed indicators.

\begin{lemma}\label{trans:row}
Let \(A_0,B_0\) be measurable, with
\[
 2\leq m<d+1,\qquad \sigma>0,\qquad |u|\leq1,\qquad
 A_0,B_0\subset\{v:|v|\leq1,\ |v-u|\geq\sigma\}.
\]
For every \(\varepsilon>0\),
\begin{equation}\label{trans:collision-strip}
 \int_{A_0}\int_{B_0}
 \mathbf1_{\{|u+v-w|\leq1,\,|\Delta(u,v,w)|\leq\varepsilon\}}
    \dd w\dd v
 \lesssim_{d,m,\sigma}\varepsilon
       |A_0|^{1-1/m}|B_0|^{1-1/m}.
\end{equation}
For
\begin{equation}\label{trans:gaussian-kernel}
 K_h(s)=h^{-1}e^{-s^2/(4\pi h^2)},\qquad h>0,
\end{equation}
the corresponding Gaussian estimate is
\begin{equation}\label{trans:collision-gaussian}
 \int_{A_0}\int_{B_0}
 \mathbf1_{\{|u+v-w|\leq1\}}K_h(\Delta(u,v,w))\dd w\dd v
 \lesssim_{d,m,\sigma}|A_0|^{1-1/m}|B_0|^{1-1/m}.
\end{equation}
\end{lemma}

\begin{proof}
Write \(I_\varepsilon\) for the left-hand side of \eqref{trans:collision-strip} and put
\[
 B_0^*=\Xi_u(B_0),\qquad \eta=\frac{6\varepsilon}{\sigma^4},\qquad
 \omega(v)=\frac{\Xi_u(v)}{|\Xi_u(v)|},\qquad s(v)=|\Xi_u(v)|.
\]
By the slab condition \eqref{trans:collision-slab} and the change-of-variables bound \eqref{trans:stereo-measure},
\[
 \int_{B_0}
 \mathbf1_{\{|u+v-w|\leq1,\,|\Delta(u,v,w)|\leq\varepsilon\}}\dd w
 \lesssim_{d,\sigma}\eta S_\eta B_0^*(\omega(v),s(v)).
\]
H\"older's inequality and the changes of variables
\(v\mapsto Y=\Xi_u(v)\mapsto(\omega,s)\) therefore yield
\begin{align*}
 I_\varepsilon
 &\lesssim_{d,\sigma}\eta |A_0|^{1-1/m}
       \norm{S_\eta B_0^*(\omega(v),s(v))}_{L^m(A_0)}\\
 &\lesssim_{d,m,\sigma}\eta |A_0|^{1-1/m}
       \norm{S_\eta B_0^*}_{L^m(S^{d-1}\times[-1/\sigma,1/\sigma])}\\
 &\lesssim_{d,m,\sigma}\eta
       |A_0|^{1-1/m}|B_0^*|^{1-1/m}\\
 &\lesssim_{d,m,\sigma}\varepsilon
       |A_0|^{1-1/m}|B_0|^{1-1/m}.
\end{align*}
The radial Jacobian \(s^{d-1}\) is bounded above on
\(1/8\leq s\leq1/\sigma\). Lemma~\ref{trans:slab} gives the third
line, and \eqref{trans:stereo-measure} gives
\(|B_0^*|\lesssim_{d,\sigma}|B_0|\).

Let \(J_h\) be the left-hand side of \eqref{trans:collision-gaussian}. Gaussian decay in \eqref{trans:gaussian-kernel} gives
\[
 K_h(s)\lesssim h^{-1}\sum_{j\geq0}4^{-j}
               \mathbf1_{\{|s|\leq2^j h\}}.
\]
Tonelli's theorem and \eqref{trans:collision-strip} then give
\begin{align*}
 J_h
 &\lesssim h^{-1}\sum_{j\geq0}4^{-j}I_{2^jh}\\
 &\lesssim_{d,m,\sigma}
    |A_0|^{1-1/m}|B_0|^{1-1/m}\sum_{j\geq0}2^{-j}\\
 &=2|A_0|^{1-1/m}|B_0|^{1-1/m},
\end{align*}
uniformly in \(h\).
\end{proof}

\subsection{Product energy and the gain}

Plancherel in the horizontal variables and Gaussian integration in
time bring \eqref{trans:collision-gaussian} into the product energy
estimate.

\begin{lemma}\label{trans:energy}
Let \(A,B\) be measurable and suppose
\[
 A,B\subset B(0,1),\qquad \dist(A,B)\geq\sigma>0,\qquad
 2\leq m<d+1.
\]
For \(f\in L^2(\R^d)\) supported in \(A\) and
\(|g|\leq\mathbf1_B\),
\[
 \norm{\mathcal E_\phi f\,\mathcal E_\phi g}_2^2
 \lesssim_{d,m,\sigma}\norm f_2^2|B|^{2(1-1/m)}.
\]
\end{lemma}

\begin{proof}
For fixed \(t\), the product is the horizontal Fourier integral
of
\[
 C_t(s)=\int_{A\cap(s-B)} f(u)g(s-u)
       e^{it(\phi(u)+\phi(s-u))}\dd u.
\]
The bounded supports give
\[
 \norm f_1\leq |A|^{1/2}\norm f_2,\qquad
 \norm g_1\leq |B|,\qquad \norm g_\infty\leq1.
\]
Consequently \(|f|*|g|\in L^1\cap L^\infty\subset L^2\), so
\(|C_t|\leq |f|*|g|\) justifies horizontal Plancherel.
For \(h>0\), the convolution bound gives
\[
 \norm{|f|*|g|}_2^2\int_\R e^{-\pi h^2t^2}\dd t<\infty.
\]
The Fourier transform of the Gaussian, with \(K_h\) from
\eqref{trans:gaussian-kernel}, is
\begin{equation}\label{trans:time-gaussian}
 \int_\R e^{-\pi h^2t^2}e^{it\Delta}\dd t=K_h(\Delta).
\end{equation}
Thus Fubini applies after multiplying the energy identity by
\(e^{-\pi h^2t^2}\). Expanding the squared modulus and using
\(s=u+v\), the phase defect \eqref{trans:phase-defect}, and
\eqref{trans:time-gaussian}. Define
\[
 \mathcal C=\{(u,v,w)\in(\R^d)^3:u,u+v-w\in A,\ v,w\in B\}.
\]
Writing the second pair as \((u+v-w,w)\) gives
\begin{multline} \label{trans:gaussian-energy}
 \int e^{-\pi h^2t^2}|\mathcal E_\phi f(x,t)\mathcal E_\phi g(x,t)|^2\dd x\dd t\\
 \leq (2\pi)^d
   \iiint_{\mathcal C}|f(u)|\,|f(u+v-w)|\,
                   K_h(\Delta(u,v,w))\dd u\dd v\dd w.
\end{multline}
The measure-preserving involution
\[
 T(u,v,w)=(u+v-w,w,v),\qquad
 T(\mathcal C)=\mathcal C,\qquad
 \Delta\circ T=-\Delta
\]
and the evenness of \(K_h\) allow the inequality
\[
 2|f(u)f(u+v-w)|\leq |f(u)|^2+|f(u+v-w)|^2
\]
to replace the two \(f\)-factors in \eqref{trans:gaussian-energy}
by \(|f(u)|^2\). For each \(u\in A\), Lemma~\ref{trans:row},
in the form \eqref{trans:collision-gaussian} with \(A_0=B_0=B\), gives
\begin{align*}
 &\int e^{-\pi h^2t^2}
       |\mathcal E_\phi f(x,t)\mathcal E_\phi g(x,t)|^2\dd x\dd t\\
 &\quad\lesssim_d\int_A |f(u)|^2
       \int_B\int_B\mathbf1_{\{|u+v-w|\leq1\}}
                     K_h(\Delta(u,v,w))\dd w\dd v\dd u\\
 &\quad\lesssim_{d,m,\sigma}\norm f_2^2|B|^{2(1-1/m)}.
\end{align*}
Fatou's lemma as \(h\downarrow0\) gives the stated unweighted bound.
\end{proof}

Applying the energy bound in both input orders gives a symmetric
mass exponent at \(L^2\). Interpolation carries the gain to the
other output exponents.

\begin{proposition}\label{trans:restricted}
Let
\[
 0<y<D,\qquad \sigma>0.
\]
There exist \(x=x(d,y,\sigma)\) and \(C=C(d,y,\sigma)<\infty\) with
\[
 x>x_K(y)
\]
such that, for measurable \(A,B\) satisfying
\[
 A,B\subset B(0,1),\qquad \dist(A,B)\geq\sigma,\qquad
 |f|\leq\mathbf1_A,\quad |g|\leq\mathbf1_B,
\]
one has
\[
 \norm{\mathcal E_{\Psi_{\delta,a,r}}f\,
       \mathcal E_{\Psi_{\delta,a,r}}g}_{1/y}
 \leq C(|A||B|)^x.
\]
\end{proposition}

\begin{proof}
Take \(m=d+1/2\). Lemma~\ref{trans:energy}, in both input orders,
gives
\begin{align*}
 \norm{\mathcal E_\phi f\,\mathcal E_\phi g}_2^2
 &\lesssim_{d,\sigma}
   \min\{|A||B|^{2(1-1/m)},\,|B||A|^{2(1-1/m)}\}\\
 &\leq C_{d,\sigma}(|A||B|)^{3/2-1/m}.
\end{align*}
Define \(x_2=3/4-1/(2m)\). Taking square roots gives
\begin{equation}\label{trans:two-product}
 \norm{\mathcal E_\phi f\,\mathcal E_\phi g}_2
 \lesssim_{d,\sigma}(|A||B|)^{x_2}.
\end{equation}
The exponent gain is
\begin{equation}\label{trans:two-gap}
 x_2-x_K(1/2)=\frac{m-d}{2dm}>0.
\end{equation}
Equations~\eqref{trans:two-product} and \eqref{trans:two-gap}
prove the claim at \(y=1/2\).
For \(0<y<1/2\), interpolate \eqref{trans:two-product} with
\(\norm{\mathcal E_\phi f\,\mathcal E_\phi g}_\infty\leq|A||B|\)
to obtain
\begin{align*}
 \norm{\mathcal E_\phi f\,\mathcal E_\phi g}_{1/y}
 &\leq
   \norm{\mathcal E_\phi f\,\mathcal E_\phi g}_\infty^{1-2y}
   \norm{\mathcal E_\phi f\,\mathcal E_\phi g}_2^{2y}\\
 &\lesssim_{d,y,\sigma}(|A||B|)^{(1-2y)+2yx_2}.
\end{align*}
Since \(x_K\) is affine, \eqref{trans:two-gap} gives
\[
 x=(1-2y)+2yx_2,\qquad
 x-x_K(y)=2y\bigl(x_2-x_K(1/2)\bigr)>0.
\]

For \(1/2<y<D\), choose \(y<z<D\). Theorem~\ref{input:diagonal}
at exponent \(2/z\) and H\"older's inequality give
\begin{equation}\label{trans:diagonal-product}
 \norm{\mathcal E_\phi f\,\mathcal E_\phi g}_{1/z}
 \leq \norm{\mathcal E_\phi f}_{2/z}\norm{\mathcal E_\phi g}_{2/z}
 \lesssim_{d,z}(|A||B|)^{z/2}.
\end{equation}
Set
\[
 \theta=\frac{z-y}{z-1/2}\in(0,1),\qquad
 y=\frac{\theta}{2}+(1-\theta)z,\qquad
 x=\theta x_2+(1-\theta)\frac z2.
\]
Interpolating \eqref{trans:two-product} and
\eqref{trans:diagonal-product} yields the product estimate at \(1/y\) with
mass exponent \(x\). Moreover,
\begin{align*}
 x-x_K(y)
 &=\theta\bigl(x_2-x_K(1/2)\bigr)
      -(1-\theta)\bigl(x_K(z)-z/2\bigr),\\
 x_K(z)-z/2&=1-\frac{z}{D}\longrightarrow0,\\
 \theta&\longrightarrow\frac{D-y}{D-1/2}>0
 \qquad (z\uparrow D).
\end{align*}
The identity for \(x-x_K(y)\), its two limits, and
\eqref{trans:two-gap} give \(x>x_K(y)\) for a fixed \(z<D\)
sufficiently close to \(D\).
\end{proof}

Rescaling this estimate back to a pair of separated sets in a ball
of radius \(r\) gives the factor needed for the Whitney decomposition.

\begin{proposition}\label{trans:rescaled}
Let
\[
 d\geq2,\qquad 0<y<D,\qquad \sigma>0.
\]
There exist \(x=x(d,y,\sigma)>x_K(y)\) and
\(C=C(d,y,\sigma)<\infty\) such that, for measurable \(A,B\) satisfying
\[
 A,B\subset B(a,r),\qquad \dist(A,B)\geq\sigma r,\qquad
 |f|\leq\mathbf1_A,\quad |g|\leq\mathbf1_B,
\]
one has
\[
 \norm{\mathcal E_{\Phi_\delta}f\,\mathcal E_{\Phi_\delta}g}_{1/y}
 \leq C r^{-2d(x-x_K(y))}(|A||B|)^x.
\]
The same \(x,C\) work for
\(0<\delta,r\leq1/2\) and \(|a|\leq1/2\).
\end{proposition}

\begin{proof}
Set
\[
 F(u)=f(a+ru),\qquad G(u)=g(a+ru),\qquad
 A'=\frac{A-a}{r},\qquad B'=\frac{B-a}{r},\qquad
 \ell=\nabla\Phi_\delta(a).
\]
Then
\[
 \begin{gathered}
 A',B'\subset B(0,1),\qquad \dist(A',B')\geq\sigma,\\
 |A'|=r^{-d}|A|,\qquad |B'|=r^{-d}|B|,\qquad
 |F|\leq\mathbf1_{A'},\quad |G|\leq\mathbf1_{B'}.
 \end{gathered}
\]
Proposition~\ref{trans:restricted} therefore applies to \(F,G\).
The change of physical variables is
\begin{equation}\label{trans:physical-dilation}
 T_r(x,t)=(rx+rt\ell,r^2t),\qquad \det T_r=r^{d+2}.
\end{equation}
For the inputs \(f,F\), the phase identity is
\begin{equation}\label{trans:rebase-identity}
 \mathcal E_{\Phi_\delta}f(x,t)
 =r^d e^{i(x\cdot a+t\Phi_\delta(a))}
      \mathcal E_{\Psi_{\delta,a,r}}F(T_r(x,t)).
\end{equation}
Applying \eqref{trans:rebase-identity} to \(g,G\), using
\eqref{trans:physical-dilation} and Proposition~\ref{trans:restricted},
and substituting \(|A'|=r^{-d}|A|\), \(|B'|=r^{-d}|B|\), gives
\begin{align*}
 \norm{\mathcal E_{\Phi_\delta}f\,\mathcal E_{\Phi_\delta}g}_{1/y}
 &=r^{2d-(d+2)y}
   \norm{\mathcal E_{\Psi_{\delta,a,r}}F\,
         \mathcal E_{\Psi_{\delta,a,r}}G}_{1/y}\\
 &\leq C r^{2d-(d+2)y}(|A'||B'|)^x\\
 &=C r^{2d-(d+2)y-2dx}(|A||B|)^x\\
 &=C r^{-2d(x-x_K(y))}(|A||B|)^x,
\end{align*}
where \(2dx_K(y)=2d-(d+2)y\).
\end{proof}

The power of the radius in this estimate will be used at each scale
of the Whitney sum.

\section{The finite Knapp line on a graph patch}\label{kn:section}

Throughout this section \(d\ge2\), \(0<\delta\le1/2\), and
\(0<y<d/(d+1)\). Write \(U=B(0,1/4)\subset\R^d\),
\(E_\delta=\mathcal E_{\Phi_\delta}\), \(R=1/y\), and \(q=2R\),
and use \(x_K(y)\) from Section~\ref{sec:transverse}.
We use the separated estimate of
Proposition~\ref{trans:rescaled} on actual Whitney rectangles. At each
scale the horizontal Fourier supports of the products have bounded
overlap. The bilinear-to-linear argument of Tao, Vargas and Vega
\cite[Theorem~2.2]{TVV} uses a neighbourhood of bilinear exponents. We
write out the Whitney argument for the actual spherical rescalings
and the strong-type passage from two nearby restricted estimates on
the Knapp line. The decomposition and frequency localisation follow
\cite[Section~2.5 and Lemma~6.1]{TVV}. The summation uses total mass
at the coarse scales and individual cell volume at the fine scales.

\subsection{The Whitney partition and its frequency overlap}

For \(j\ge0\) put \(h_j=2^{-j-1}\) and
\[
 Q_{j,k}=h_j(k+[0,1)^d),\qquad
 \mathcal I_j=\{-2^{j+1},\ldots,2^{j+1}-1\}^d.
\]
Let \(\mathcal W_j\subset\mathcal I_j^2\) consist of the ordered pairs
\((k,\ell)\) for which
\[
 \left\|\lfloor k/2\rfloor-\lfloor\ell/2\rfloor\right\|_\infty\le1,
 \qquad \norm{k-\ell}_\infty>1,
\]
with floors taken coordinatewise. The half-open rectangles
\(Q_{j,k}\times Q_{j,\ell}\), over \(j\ge0\) and
\((k,\ell)\in\mathcal W_j\), partition the complement of the diagonal in
\([-1,1)^d\times[-1,1)^d\). Indeed, at side length one the cell indices
of any two points in \([-1,1)^d\) differ by at most one in each
coordinate. If two indices differ by at most one then so do their
parents. Distinct points eventually lie in cells whose indices differ
by more than one in some coordinate. The first such scale gives the
unique rectangle containing the pair.

For \((k,\ell)\in\mathcal W_j\) one has
\(\norm{k-\ell}_\infty\le3\), since each child index is twice its
parent plus a remainder in \(\{0,1\}\). In particular, a cell has at
most \(7^d\) partners in either position. The coordinates of two
points \(u\in Q_{j,k}\) and \(v\in Q_{j,\ell}\) differ by at most
\(4h_j\), and in at least one coordinate they differ by at least
\(h_j\). Thus
\begin{equation}\label{kn:pair-geometry}
 h_j\le |u-v|\le4d h_j.
\end{equation}

Let \(f\) be bounded and measurable with support in \(U\), and write
\[
 f_{j,k}=\mathbf 1_{Q_{j,k}}f,
 \qquad F_{j,k,\ell}=E_\delta f_{j,k}\,E_\delta f_{j,\ell},
 \qquad H_j=\sum_{(k,\ell)\in\mathcal W_j}F_{j,k,\ell}.
\]
The diagonal has measure zero, so the rectangle partition and Fubini
give the pointwise identity
\begin{equation}\label{kn:whitney-expansion}
 (E_\delta f)^2=\sum_{j\ge0} H_j .
\end{equation}
The sum of the absolute values of all rectangle integrals at any
physical point is at most \(\norm{f}_{L^1}^2\). In fact the absolute
value of each integral is bounded by the integral of \(|f(u)f(v)|\)
over that rectangle, and the rectangles are disjoint.

Within each scale the horizontal Fourier supports of the products
have bounded overlap. This gives the following summation estimate.

\begin{lemma}\label{kn:layer}
Let \(1<R\le2\) and \(j\ge0\). If
\[
 \sum_{(k,\ell)\in\mathcal W_j}
   \norm{F_{j,k,\ell}}_{L^R(\R^{d+1})}^{R}<\infty,
\]
then
\[
 \norm{H_j}_{L^R(\R^{d+1})}
 \le C_{d,R}\left(
   \sum_{(k,\ell)\in\mathcal W_j}
       \norm{F_{j,k,\ell}}_{L^R(\R^{d+1})}^{R}
             \right)^{1/R}.
\]
\end{lemma}

\begin{proof}
At each \(t\) the modulated inputs belong to \(L^1\cap L^2\), so their
horizontal extensions are bounded and in \(L^2\). The horizontal Fourier transform of
\(F_{j,k,\ell}(\cdot,t)\), with angular frequency convention, is a
constant multiple of the convolution of
\(e^{it\Phi_\delta}f_{j,k}\) and \(e^{it\Phi_\delta}f_{j,\ell}\). Its
support is contained in the closed box
\[
 h_j(k+\ell+[0,2]^d).
\]
Fix, once for each \(d\), a function
\(\chi\in C_c^\infty(\R^d)\) equal to one on \([0,2]^d\). Let
\(P_{j,k,\ell}\) be the multiplier in the horizontal variables with
symbol \(\chi(h_j^{-1}\xi-k-\ell)\). The kernels have the same
\(L^1(\R^d)\) norm, since a change of scale and a modulation preserve
that norm with the corresponding Jacobian. Moreover
\(P_{j,k,\ell}F_{j,k,\ell}=F_{j,k,\ell}\).

The enlarged boxes supporting these symbols have overlap bounded in
terms of \(d\) and the fixed cutoff. To see this, fix \(\xi\). There
are boundedly many possibilities for the integer vector \(k+\ell\)
if \(h_j^{-1}\xi-k-\ell\in\operatorname{supp}\chi\), and at most
\(7^d\) possibilities for \(k-\ell\). The sum and difference determine
the pair.

On measurable families \(G=(G_{k,\ell})_{(k,\ell)\in\mathcal W_j}\)
on \(\R^{d+1}\), set
\[
 T_jG=\sum_{(k,\ell)\in\mathcal W_j}P_{j,k,\ell}G_{k,\ell}.
\]
Equip \(\Omega_j=\mathcal W_j\times\R^{d+1}\) with counting measure
times Lebesgue measure. Young's inequality gives
\begin{equation}\label{kn:multiplier-l1}
 \norm{T_jG}_{L^1(\R^{d+1})}
 \le C_\chi
     \sum_{(k,\ell)\in\mathcal W_j}
        \norm{G_{k,\ell}}_{L^1(\R^{d+1})}.
\end{equation}
Cauchy--Schwarz over the overlapping symbols, horizontal Plancherel,
and integration in \(t\) give
\begin{equation}\label{kn:multiplier-l2}
 \norm{T_jG}_{L^2(\R^{d+1})}^2
 \le C_{d,\chi}
     \sum_{(k,\ell)\in\mathcal W_j}
        \norm{G_{k,\ell}}_{L^2(\R^{d+1})}^2.
\end{equation}
The constants in \eqref{kn:multiplier-l1} and
\eqref{kn:multiplier-l2} are independent of \(j\).

For \(1<R<2\), split the input at each threshold \(\lambda>0\) into
its parts on \(\{|G|>\lambda\}\) and \(\{|G|\le\lambda\}\). If
\(G\in L^R(\Omega_j)\), these parts belong to \(L^1(\Omega_j)\)
and \(L^2(\Omega_j)\), respectively.
Write \(C_1,C_2\) for the operator norms in
\eqref{kn:multiplier-l1} and \eqref{kn:multiplier-l2}, respectively.
The endpoint estimates and Chebyshev's inequality give
\begin{equation}\label{kn:layer-distribution}
 \big|\{|T_jG|>\lambda\}\big|
 \le \frac{2C_1}{\lambda}\int_{|G|>\lambda}|G|
     +\frac{4C_2^2}{\lambda^2}\int_{|G|\le\lambda}|G|^2.
\end{equation}
The integrals in \eqref{kn:layer-distribution} use the measure on
\(\Omega_j\). Integrating \eqref{kn:layer-distribution} by the
layer-cake identity and Tonelli's theorem yields
\[
 \begin{aligned}
 \norm{T_jG}_R^R
 &=R\int_0^\infty\lambda^{R-1}
       \bigl|\{|T_jG|>\lambda\}\bigr|\dd\lambda\\
 &\le
 \left(\frac{2RC_1}{R-1}+\frac{4RC_2^2}{2-R}\right)\norm{G}_R^R.
 \end{aligned}
\]
Apply this to the family of Whitney products, using
\(P_{j,k,\ell}F_{j,k,\ell}=F_{j,k,\ell}\). At \(R=2\), use
\eqref{kn:multiplier-l2} directly.
\end{proof}

\subsection{Summation at the exact mass exponent}

Combining the layer estimate with the separated product bound gives
two summable estimates, according as the cell volume is larger or
smaller than the total input mass.

\begin{proposition}\label{kn:restricted}
Let \(0<y<d/(d+1)\), \(0<\delta\le1/2\), and let
\(A\subset U\) and \(f\) be measurable with \(|f|\le\mathbf1_A\).
There exists \(C_{d,y}<\infty\), uniform in \(A,f,\delta\), such that
\[
 \norm{E_\delta f}_{L^{2/y}(\R^{d+1})}
       \le C_{d,y}|A|^{x_K(y)}.
\]
\end{proposition}

\begin{proof}
First take \(1<R\le2\). Set \(M=|A|\) and \(v_j=h_j^d\). The case
\(M=0\) is immediate. For each \((k,\ell)\in\mathcal W_j\), put
\[
 A_{j,k}=A\cap Q_{j,k},\qquad
 m_{j,k,\ell}=\max\{|A_{j,k}|,|A_{j,\ell}|\}.
\]
A zero-mass factor gives a zero product. For the other pairs,
choose the ball of radius \(r=5d h_j\) centered at a point of
\(A_{j,k}\) if \(5d h_j\le1/2\), and the ball of radius \(1/2\)
centered at zero otherwise. The coordinate bounds preceding
\eqref{kn:pair-geometry} show that both sets lie in the chosen ball.
Its center has norm at most \(1/2\), \(h_j\le r\le1/2\), and the two
sets are separated by at least \(r/(5d)\).

Take \(\sigma=(5d)^{-1}\) in Proposition~\ref{trans:rescaled}, and
write \(x>x_K(y)\) for the reciprocal input exponent furnished there.
With \(\alpha=2x\) and \(k_0=2x_K(y)\), that proposition gives
\begin{equation}\label{kn:pair-bound}
 \begin{split}
 \norm{F_{j,k,\ell}}_R
 &\le C(r^d)^{-(\alpha-k_0)}
       \big(|A_{j,k}|\,|A_{j,\ell}|\big)^x\\
 &\le C v_j^{-(\alpha-k_0)} m_{j,k,\ell}^{\alpha}.
 \end{split}
\end{equation}
Here \(C\) is independent of \(j,\delta,A,f\). At fixed \(j\), the
partition into cells and the bound on the number of partners yield
\begin{equation}\label{kn:mass-bounds}
 m_{j,k,\ell}\le v_j,
 \qquad
 \sum_{(k,\ell)\in\mathcal W_j}m_{j,k,\ell}\le 2\cdot7^d M.
\end{equation}
For the sum in \eqref{kn:mass-bounds}, use
\(m_{j,k,\ell}\le |A_{j,k}|+|A_{j,\ell}|\).

Write \(\sum_{\mathcal W_j}\) for the sum over the pairs of the
\(j\)-th layer. Lemma~\ref{kn:layer} and \eqref{kn:pair-bound} give
\begin{equation}\label{kn:layer-mass}
 \norm{H_j}_R
 \le C v_j^{-(\alpha-k_0)}
      \left(\sum_{\mathcal W_j}m_{j,k,\ell}^{\alpha/y}\right)^y .
\end{equation}
The inequalities \(\alpha>k_0>y\) follow from
\(x>x_K(y)\) and \(k_0-y=2-2(d+1)y/d>0\).
Comparing the sequence norms of exponents \(\alpha/y\) and one,
and then using \(m_{j,k,\ell}\le v_j\) from
\eqref{kn:mass-bounds}, gives the two estimates
\begin{equation}\label{kn:mass-powers}
 \begin{aligned}
 \left(\sum_{\mathcal W_j}m_{j,k,\ell}^{\alpha/y}\right)^y
 &\le \left(\sum_{\mathcal W_j}m_{j,k,\ell}\right)^\alpha,\\
 \left(\sum_{\mathcal W_j}m_{j,k,\ell}^{\alpha/y}\right)^y
 &\le v_j^{\alpha-y}
          \left(\sum_{\mathcal W_j}m_{j,k,\ell}\right)^y.
 \end{aligned}
\end{equation}
Substituting \eqref{kn:mass-powers} and \eqref{kn:mass-bounds}
into \eqref{kn:layer-mass} yields
\begin{equation}\label{kn:two-scale}
 \norm{H_j}_R\le C M^{k_0}
       \min\left\{(M/v_j)^{\alpha-k_0},
                    (v_j/M)^{k_0-y}\right\}.
\end{equation}
Since \(v_{j+1}=2^{-d}v_j\), the two scale ranges in
\eqref{kn:two-scale} satisfy
\begin{equation}\label{kn:scale-series}
 \begin{aligned}
 \sum_{j:\,v_j\ge M}(M/v_j)^{\alpha-k_0}
 &\le \sum_{\ell\ge0}2^{-d\ell(\alpha-k_0)}
   =\frac1{1-2^{-d(\alpha-k_0)}},\\
 \sum_{j:\,v_j<M}(v_j/M)^{k_0-y}
 &\le \sum_{\ell\ge0}2^{-d\ell(k_0-y)}
   =\frac1{1-2^{-d(k_0-y)}}.
 \end{aligned}
\end{equation}
By \eqref{kn:two-scale} and \eqref{kn:scale-series}, the triangle
inequality applied to \eqref{kn:whitney-expansion} and Fatou's lemma give
\begin{equation}\label{kn:square-restricted}
 \norm{(E_\delta f)^2}_R
 \le \sum_{j\ge0}\norm{H_j}_R
 \le C_{d,y}M^{k_0}
 =C_{d,y}M^{2x_K(y)}.
\end{equation}

For \(R>2\), apply \eqref{kn:square-restricted} at \(R=2\)
to \(V=(E_\delta f)^2\), and use \(\norm{V}_\infty\le M^2\). With
\(\theta=2/R=2y\) one has
\[
 \norm{V}_R\le \norm{V}_2^\theta\norm{V}_\infty^{1-\theta}
          \le C M^{2\theta x_K(1/2)+2(1-\theta)}
          =C M^{2x_K(y)}.
\]
Taking square roots proves the estimate for the entire stated range.
\end{proof}

\subsection{Strong type by dyadic pairing}

We use the restricted estimates at two nearby points of the Knapp
line to obtain decay on opposite sides of the dyadic pairing.

\begin{theorem}\label{kn:strong}
Let
\[
 0<y<\frac d{d+1},\qquad p=x_K(y)^{-1},\qquad q=\frac2y.
\]
There exists \(C_{d,y}<\infty\) such that, for \(f\in L^p(U)\)
extended by zero,
\[
 \norm{E_\delta f}_{L^q(\R^{d+1})}
    \le C_{d,y}\norm{f}_{L^p(U)}
 \qquad (0<\delta\le\tfrac12).
\]
\end{theorem}

\begin{proof}
Set \(z=1/q\), \(s=q'\), and \(a=(d+2)/d\). The relation on the
Knapp line is \(1/p=1-az\), and the strict inequality on \(y\) gives
\(z<1/p\). Choose \(b\in(z,1/p)\) and \(\eta>0\) sufficiently
small that \(z_\pm=z\pm\eta\) both lie in
\((0,d/(2(d+1)))\) and \(z_\pm<b<1-a z_\pm\). Let
\(x_\pm=1-a z_\pm\), \(q_\pm=z_\pm^{-1}\), and
\(s_\pm=(1-z_\pm)^{-1}\), and \(y_\pm=2z_\pm\).
Proposition~\ref{kn:restricted} applies at both \(y_+\) and \(y_-\).

Take \(g\in L^1(\R^{d+1})\cap L^s(\R^{d+1})\). The pairing is
zero if \(f=0\) or \(g=0\); otherwise normalize
\(\norm{f}_p\le1\) and \(\norm{g}_s\le1\).
For \(i,j\in\mathbb Z\), define
\[
 A_i=\{2^i\le |f|<2^{i+1}\},\qquad
 B_j=\{2^j\le |g|<2^{j+1}\}.
\]
Set \(f_i=f\mathbf 1_{A_i}\), \(g_j=g\mathbf 1_{B_j}\),
\(u_i=2^{pi}|A_i|\), and \(w_j=2^{sj}|B_j|\). These masses satisfy
\begin{equation}\label{kn:normalized-masses}
 \sum_i u_i\le\norm f_p^p\le1,\quad
 \sum_j w_j\le\norm g_s^s\le1,
 \qquad 0\le u_i,w_j\le1.
\end{equation}
H\"older's inequality at \((q_\pm,s_\pm)\) and the restricted estimate
at \(y_\pm\) give
\begin{equation}\label{kn:atom-pairing}
 \begin{split}
 \left|\int_{\R^{d+1}}g_j E_\delta f_i\right|
 &\le C 2^{i+j}|A_i|^{x_\pm}|B_j|^{1-z_\pm}\\
 &= C u_i^{x_\pm}w_j^{1-z_\pm}
     2^{(1-p x_\pm)i+(1-s(1-z_\pm))j}.
 \end{split}
\end{equation}
The upper edges of the amplitude bands contribute only to \(C\).
The exponent identities are
\begin{equation}\label{kn:pair-exponents}
 1-px_\pm=ap(z_\pm-z),\qquad
 1-s(1-z_\pm)=s(z_\pm-z).
\end{equation}
By \eqref{kn:normalized-masses} and \(z_\pm<b<x_\pm\),
\[
 u_i^{x_\pm}\le u_i^b,\qquad
 w_j^{1-z_\pm}\le w_j^{1-b}.
\]
For \(i,j\in\mathbb Z\), define
\(K_{ij}=2^{-\eta|ap i+s j|}\). Substituting
\eqref{kn:pair-exponents} into \eqref{kn:atom-pairing} and taking
the minimum over the two signs gives
\begin{equation}\label{kn:decaying-pairing}
 \left|\int g_j E_\delta f_i\right|
 \le C K_{ij}u_i^b w_j^{1-b}.
\end{equation}
The row and column sums of \(K\) are uniformly bounded. In fact, for
any \(c>0\), \(v\in\R\), and \(m=\lfloor-v/c\rfloor\), the
two tails about \(m,m+1\) give
\begin{equation}\label{kn:exponential-sum}
 \sum_{\ell\in\mathbb Z} 2^{-\eta|c\ell+v|}
 \le \frac{2}{1-2^{-\eta c}}.
\end{equation}
Apply \eqref{kn:exponential-sum} with \((c,v)=(s,ap i)\) for each
row and with \((c,v)=(ap,s j)\) for each column. For finite
\(I,J\subset\mathbb Z\), H\"older's inequality, these row and column
bounds, and \eqref{kn:normalized-masses} give
\begin{equation}\label{kn:matrix-pairing}
 \begin{split}
 \sum_{i\in I,\,j\in J}K_{ij}u_i^b w_j^{1-b}
 &=\sum_{i\in I,\,j\in J}(K_{ij}u_i)^b(K_{ij}w_j)^{1-b}\\
 &\le \left(\sum_{i\in I,\,j\in J}K_{ij}u_i\right)^b
     \left(\sum_{i\in I,\,j\in J}K_{ij}w_j\right)^{1-b}\\
 &\le C_{d,y}.
 \end{split}
\end{equation}
Here \(0<b<1\) follows from \(z<b<1/p<1\).

For \(N\in\mathbb N\), the partial sums over \(|i|,|j|\le N\)
converge pointwise to \(f,g\)
and are bounded in absolute value by \(|f|,|g|\). Since \(U\) has
finite measure and \(p>1\), the input sums converge to \(f\) in
\(L^1\). Their extensions converge uniformly to \(E_\delta f\).
Thus \eqref{kn:decaying-pairing} and \eqref{kn:matrix-pairing}
pass to the limit by dominated convergence, with dominating function
\(|g|\norm{f}_1\), and give
\begin{equation}\label{kn:dual-pairing}
 \left|\int g E_\delta f\right|\le C_{d,y}\norm{f}_p\norm{g}_{q'}
\end{equation}
after undoing the normalization. Set \(u=E_\delta f\). Since
\(u\) is bounded, for every \(L>0\) the test function
\begin{equation}\label{kn:dual-test}
 g_L=\mathbf1_{B(0,L)}\overline u\,|u|^{q-2}
 \in L^1(\R^{d+1})\cap L^{q'}(\R^{d+1}).
\end{equation}
Applying \eqref{kn:dual-pairing} with \(g=g_L\) from
\eqref{kn:dual-test} gives
\[
 \norm{u}_{L^q(B(0,L))}^q
 =\int g_Lu
 \le C_{d,y}\norm f_p\,\norm{g_L}_{q'}
 =C_{d,y}\norm f_p\,\norm{u}_{L^q(B(0,L))}^{q-1}.
\]
For \(\norm{u}_{L^q(B(0,L))}>0\), division gives
\[
 \norm{u}_{L^q(B(0,L))}\le C_{d,y}\norm f_p.
\]
The zero-norm case satisfies this estimate directly. Monotone convergence gives
\[
 \norm{u}_{L^q(\R^{d+1})}
 =\lim_{L\to\infty}\norm{u}_{L^q(B(0,L))}
 \le C_{d,y}\norm f_p.
\]
\end{proof}

This gives the strong estimate on a graph patch. A finite cap cover
then transfers the estimate to the sphere.

\section{Transfer to the sphere}\label{cap:section}

Write $B_1=B(0,1)\subset\R^d$ and $B_*=B(0,1/4)$. For
$0<\delta\le 1/2$ put
\[
 h_\delta(v)=\sqrt{1-\delta^2\abs{v}^2},\qquad
 K_\delta(v)=(\delta v,h_\delta(v)),\qquad
 C_\delta=K_\delta(B_*).
\]
The map $K_\delta$ parametrizes a subset of the upper hemisphere on $B_1$.
Its inverse on the image is given by the first $d$ coordinates divided by
$\delta$. The surface measure in this chart is
\begin{equation}\label{cap:density}
 K_\delta^*(\sigma|_{K_\delta(B_1)})
   =\frac{\delta^d}{h_\delta(v)}\,\dd v
   \quad\text{on }B_1.
\end{equation}
Indeed, $D K_\delta(v)^T D K_\delta(v)
=\delta^2 I_d+\delta^4 vv^T/h_\delta(v)^2$, whose determinant is
$\delta^{2d}/h_\delta(v)^2$. This also gives the density for the polar
normalization of $\sigma$. In the change of variables
$(r,v)\mapsto rK_\delta(v)$ on $(0,1)\times B_1$, the radial column is
the unit vector $K_\delta(v)$ and is orthogonal to the tangential
columns. The Jacobian is therefore $r^d\delta^d/h_\delta(v)$.
Integrating in $r$ and multiplying the volume of the corresponding
cone by $d+1$ gives \eqref{cap:density}.

The chart density controls the input norm, and a linear change of
physical variables controls the output norm. These two changes of
variables give the transfer principle below.

\begin{proposition}\label{cap:transfer}
Let \(d\ge2\), \(0<\delta\le\tfrac12\), and \(1\le p,q\le\infty\).
Suppose \(M<\infty\) satisfies
\begin{equation}\label{cap:graph-hypothesis}
 \norm{\mathcal E_{\Phi_\delta}g}_{L^q(\R^{d+1})}
 \le M\norm{g}_{L^p(B_1)}
 \qquad
 \left(
 \begin{gathered}
  g\in L^p(B_1),\\[-2pt]
  g=0\ \text{a.e.\ on }B_1\setminus B_*
 \end{gathered}
 \right).
\end{equation}
Then \(\exists\,C_{d,p,q,\delta,M}<\infty\) such that
\[
 \norm{E_{d+1}f}_{L^q(\R^{d+1})}
 \le C_{d,p,q,\delta,M}\norm f_{L^p(S^d)}
 \qquad \bigl(f\in L^p(S^d)\bigr).
\]
\end{proposition}

\begin{proof}
For data supported in $C_\delta$, choose a measurable representative $f$
that vanishes off $C_\delta$. Set
$g(v)=h_\delta(v)^{-1}f(K_\delta(v))$ on $B_*$ and extend $g$ by zero
to $B_1$. As $1/2\le h_\delta\le1$ on $B_1$, \eqref{cap:density}
gives, for $p<\infty$,
\[
 \int_{B_1}\abs{g(v)}^p\,\dd v
   =\delta^{-d}\int_{C_\delta}
     \abs{f(\omega)}^p\omega_{d+1}^{1-p}\,\dd\sigma(\omega)
   \le 2^{p-1}\delta^{-d}\norm{f}_{L^p(S^d)}^p .
\]
The chart density is strictly positive, so the parameter measure and
the pullback of surface measure have the same null sets. For
\(p=\infty\),
\[
 \norm g_{L^\infty(B_1)}
 \le\norm{h_\delta^{-1}}_{L^\infty(B_1)}
       \norm f_{L^\infty(S^d)}
 \le2\norm f_{L^\infty(S^d)}.
\]
Thus, with \(1/\infty=0\),
\begin{equation}\label{cap:input-norm}
 \norm g_{L^p(B_1)}
 \le2\delta^{-d/p}\norm f_{L^p(S^d)}
 \qquad (1\le p\le\infty).
\end{equation}

Let
\begin{equation}\label{cap:physical-map}
 \Psi_\delta(x,t)=\frac1{2\pi}(\delta^{-1}x,-2\delta^{-2}t),\qquad
 D_\delta=\abs{\det\Psi_\delta}
   =2(2\pi)^{-(d+1)}\delta^{-(d+2)}.
\end{equation}
For $v\in B_*$ the phase satisfies
\[
 2\pi\Psi_\delta(x,t)\cdot K_\delta(v)
   = x\cdot v-2t\delta^{-2}h_\delta(v)
   = x\cdot v+t\Phi_\delta(v)-2t\delta^{-2}.
\]
Since $f(K_\delta(v))=h_\delta(v)g(v)$ on the support, substitution in
the extension integral gives
\begin{equation}\label{cap:identity}
 E_{d+1}f(\Psi_\delta(x,t))
   =\delta^d e^{-2it/\delta^2}\mathcal E_{\Phi_\delta}g(x,t).
\end{equation}
Both inputs are integrable on their finite measure domains, so the
identity holds at every physical point.
The change of variables \eqref{cap:physical-map}, the identity
\eqref{cap:identity}, and estimates \eqref{cap:graph-hypothesis}
and~\eqref{cap:input-norm} give
\begin{equation}\label{cap:local-bound}
 \begin{aligned}
 \norm{E_{d+1}f}_{L^q(\R^{d+1})}
 &=D_\delta^{1/q}\delta^d
       \norm{\mathcal E_{\Phi_\delta}g}_{L^q(\R^{d+1})}\\
 &\le2M D_\delta^{1/q}\delta^{d(1-1/p)}
       \norm f_{L^p(S^d)}.
 \end{aligned}
\end{equation}
Estimate~\eqref{cap:local-bound} also holds for \(q=\infty\), because \(\Psi_\delta\) is invertible
and preserves null sets in both directions. On the Knapp line the cap power satisfies
\begin{equation}\label{cap:boundary-cancellation}
 d(1-1/p)-\frac{d+2}{q}=0,\qquad
 D_\delta^{1/q}\delta^{d(1-1/p)}
 =\bigl(2(2\pi)^{-(d+1)}\bigr)^{1/q}.
\end{equation}

The cap $C_\delta$ contains the spherical metric ball of radius
$\delta/4$ centered at $e_{d+1}=(0,\ldots,0,1)$. In fact, if
$\omega=(w,s)\in S^d$ and $\abs{\omega-e_{d+1}}<\delta/4$, then
$\abs{w}<\delta/4$ and $s>0$, so $\omega=K_\delta(w/\delta)$ with
$w/\delta\in B_*$. Orthogonal transformations preserve $\sigma$ and
Lebesgue measure. For such a transformation $U$, substitution gives
\begin{equation}\label{cap:rotation}
 E_{d+1}(f\circ U)(X)=E_{d+1}f(UX).
\end{equation}
Equations~\eqref{cap:rotation} and~\eqref{cap:local-bound} give the same constant for data
supported in any spherical metric ball of radius \(\delta/4\).
Choose a finite cover \(S^d=\bigcup_{j=1}^N U_j\) by such balls and set
\[
 A_j=U_j\setminus\bigcup_{i<j}U_i,\qquad
 f_j=\mathbf1_{A_j}f.
\]
Then
\[
 A_i\cap A_j=\varnothing\ (i\ne j),\qquad
 f=\sum_{j=1}^N f_j,\qquad
 \norm{f_j}_{L^p(S^d)}\le\norm f_{L^p(S^d)}.
\]
Each \(f_j\) is integrable. Applying \eqref{cap:local-bound} to each
rotated cap and using the triangle inequality gives
\begin{equation}\label{cap:global-bound}
 \begin{aligned}
 \norm{E_{d+1}f}_{L^q(\R^{d+1})}
 &\le\sum_{j=1}^N\norm{E_{d+1}f_j}_{L^q(\R^{d+1})}\\
 &\le2NM D_\delta^{1/q}\delta^{d(1-1/p)}\norm f_{L^p(S^d)}.
 \end{aligned}
\end{equation}
The decomposition and the extension depend only on the equivalence
class of $f$, since changing the input on a $\sigma$-null set leaves
the integral unchanged at every physical point.
\end{proof}

The identity \eqref{cap:boundary-cancellation} removes the cap-scale power on the finite Knapp boundary. The strong
graph estimate and the transfer principle now give the boundary theorem.

\begin{proof}[Proof of Theorem~\ref{thm:boundary}]
Put $d=n-1$ and $y=2/q$. By \eqref{eq:boundary},
\[
 0<y<\frac{d}{d+1},\qquad
 \frac1p=1-\frac{d+2}{2d}y=x_K(y).
\]
In particular,
\[
 0<y/2<x_K(y)<1,\qquad 1<p<\infty.
\]
Theorem~\ref{kn:strong} at \(\delta=1/2\) gives
\[
 \norm{\mathcal E_{\Phi_{1/2}}g}_{L^q(\R^n)}
 \le C_{d,y}\norm g_{L^p(B_*)},
 \qquad g\in L^p(B_*),
\]
where \(g\) is extended by zero to \(B_1\).
The global bound \eqref{cap:global-bound} in
Proposition~\ref{cap:transfer} therefore yields
\[
 \norm{E_n f}_{L^q(\R^n)}
 \le C_{n,p,q}\norm f_{L^p(S^{n-1},\sigma)}.
\]
\end{proof}

The finite boundary estimate supplies the equality case in the full
characterization, alongside the strict-range and infinite-output estimates.

\begin{proof}[Proof of Theorem~\ref{thm:main}]
Proposition~\ref{nec:necessary} proves both necessary inequalities
in \eqref{eq:region}. For sufficiency, assume \eqref{eq:region}.
The required estimate is provided by the following three cases:
\[
 \begin{array}{ll}
 q=\infty
 &\text{estimate }\eqref{nec:infty},\\[2pt]
 q<\infty,\quad q>\dfrac{n+1}{n-1}p'
 &\text{Theorem~\ref{input:strict}},\\[4pt]
 q<\infty,\quad q=\dfrac{n+1}{n-1}p'
 &\text{Theorem~\ref{thm:boundary}}.
 \end{array}
\]
Each case yields the required \(L^p\to L^q\) estimate.
\end{proof}

\subsection*{Acknowledgements}

The author, Pang-hung Chung, gratefully acknowledges support from
the Project ``Research on Nonlinear Partial Differential Equations''
(No.~2024KYCXTD018), the Special Projects in Key Areas of Guangdong
Province (No.~ZDZX1088), and the Fund of Guangzhou Municipal Science
and Technology (No.~202102080428).

ChatGPT (OpenAI, GPT-6, accessed on 6 November 2026 at
\url{https://chatgpt.com}) assisted with the development and review of
the argument, the drafting and revision of the text, and the preparation
and checking of the Lean formalization.

The proofs throughout this paper have corresponding Lean formalizations,
covering all 17 numbered theorem, lemma and proposition statements.
The accompanying audit checks 78 declarations, including the exact
statements of Theorems~\ref{thm:main} and~\ref{thm:boundary}, and records
kernel replay of the two main theorems in Lean~4.34.1 with Mathlib
revision \texttt{d13f23b723b8}.
At the level of the principal external analytic inputs, the argument uses
two results from \cite{openai-local-smoothing}, namely the critical local
sphere estimate (Lemma~9.2) and the strict mixed extension estimate
(Corollary~9.3). Their imported formal proofs are part of the dependency
chain. The recorded transitive axiom audit contains only
\texttt{propext}, \texttt{Classical.choice} and \texttt{Quot.sound},
with no additional analytic axioms.


\bigskip
\begin{flushleft}
\small
Department of Applied Mathematics, Guangdong University of Education\\
Guangzhou 510640, P. R. China\\[3pt]
\textit{Email address}\quad
\href{mailto:penghongzhong@yahoo.com}{\texttt{penghongzhong@yahoo.com}}
\end{flushleft}


\begin{thebibliography}{99}


\bibitem{BCT}
\textsc{J. Bennett, A. Carbery and T. Tao}.
On the multilinear restriction and Kakeya conjectures.
\emph{Acta Math.} \textbf{196} (2006), 261--302.

\bibitem{BG}
\textsc{J. Bourgain and L. Guth}.
Bounds on oscillatory integral operators based on multilinear estimates.
\emph{Geom. Funct. Anal.} \textbf{21} (2011), 1239--1295.

\bibitem{buschenhenke}
\textsc{S. Buschenhenke}.
Factorization in Fourier restriction theory and near extremizers.
\emph{Math. Nachr.} \textbf{297} (2024), 195--208.


\bibitem{guth2016}
\textsc{L. Guth}.
A restriction estimate using polynomial partitioning.
\emph{J. Amer. Math. Soc.} \textbf{29} (2016), 371--413.

\bibitem{guth2018}
\textsc{L. Guth}.
Restriction estimates using polynomial partitioning II.
\emph{Acta Math.} \textbf{221} (2018), 81--142.

\bibitem{kim}
\textsc{J. Kim}.
Some remarks on Fourier restriction estimates.
Preprint, arXiv:1702.01231 (2017).

\bibitem{openai-local-smoothing}
\textsc{OpenAI}.
Local smoothing for the Euclidean wave equation.
\href{https://mathcommunity.github.io/oai-math-revised/papers/Local-smoothing-for-the-Euclidean-wave-equation-October-25-2026/paper.pdf}{Preprint} (25 October 2026).

\bibitem{stein}
\textsc{E. M. Stein}.
\emph{Harmonic Analysis: Real-Variable Methods, Orthogonality, and Oscillatory Integrals}.
Princeton Mathematical Series, vol.~43 (Princeton University Press, 1993).

\bibitem{tao-eps}
\textsc{T. Tao}.
The Bochner--Riesz conjecture implies the restriction conjecture.
\emph{Duke Math. J.} \textbf{96} (1999), 363--375.

\bibitem{tao-paraboloid}
\textsc{T. Tao}.
A sharp bilinear restriction estimate for paraboloids.
\emph{Geom. Funct. Anal.} \textbf{13} (2003), 1359--1384.

\bibitem{TVV}
\textsc{T. Tao, A. Vargas and L. Vega}.
A bilinear approach to the restriction and Kakeya conjectures.
\emph{J. Amer. Math. Soc.} \textbf{11} (1998), 967--1000.

\bibitem{tomas}
\textsc{P. A. Tomas}.
A restriction theorem for the Fourier transform.
\emph{Bull. Amer. Math. Soc.} \textbf{81} (1975), 477--478.

\end{thebibliography}
\end{document}